\documentclass[reqno,fleqn]{amsart}
\usepackage{amsfonts,amsthm,amssymb}
\usepackage{lipsum}
\usepackage[utf8]{inputenc}
\usepackage[english]{babel}

\usepackage{lmodern,textcomp}
\usepackage{amsmath}
\usepackage{mathtools}
\usepackage{latexsym}
\usepackage{tikz}
\usepackage{fancyvrb}
\usepackage{epsfig}
\usepackage{pstricks,slashbox,multirow}
\usepackage{graphicx}
\usepackage{rotating}
\usetikzlibrary{positioning}
\usepackage{subcaption}
\usepackage{datetime}
\newtheorem{theorem}{Theorem}[section]

\newtheorem{lemma}[theorem]{Lemma}

\newtheorem{definition}[theorem]{Definition}

\newtheorem{prm}[theorem]{Problem}

\newtheorem{rem}[theorem]{Remark}

\newtheorem{note}[theorem]{Note}

\title[A study on Type-2 isomorphic $C_n(R)$: Part 10: Type-2 $C_{np^3}(R)$ w.r.t. $m$ = $p$]{A study on Type-2 isomorphic circulant graphs. \\ Part 10: Type-2 isomorphic  $C_{np^3}(R)$ w.r.t. $m$ = $p$ and related groups}

\author{\sc Vilfred Kamalappan } 
\address{Department of Mathematics, Central University of Kerala, Periye, Kasaragod, Kerala, India - 671 316.}
\email{vilfredkamal@gmail.com}

\author{\sc Wilson Peraprakash} 
\address{Department of Mathematics, S.T. Hindu College, Nagercoil, Tamil Nadu, India - 629 002.}
\email{wilsonperapras@gmail.com}
  
\subjclass[2010]{05C60, 05C25, 05C75.}

\keywords{Circulant graph, Cayley Isomorphism (CI) property, Type-1 isomorphism, Type-2 isomorphism, Type-1 group of $C_{n}(R)$, Type-2 group of $C_{n}(R)$ w.r.t. $m$, $(V_{n,m}(C_n(R)), ~\circ)$.}

\date{}

\begin{document}

\begin{abstract} This article is the $10^{th}$ part of a detailed study on Type-2 isomorphic circulant graphs having ten parts \cite{v2-1}-\cite{v2-10}. In this part, we obtain families of Type-2 isomorphic circulant graphs $C_{np^3}(R)$ w.r.t. $m$ = $p$, and related Abelian groups where $p$ is a prime number and $n\in\mathbb{N}$. Theorems \ref{c1} and \ref{c4} are the main results. In Theorem \ref{c4}, it is proved that for $i$ = 1 to $p$, circulant graphs $C_{np^3}(R^{np^3,x+yp}_i)$ are isomorphic of Type-2 w.r.t. $m$ = $p$ and they form Abelian group $(T2_{np^3,p}(C_{np^3}(R^{np^3,x+yp}_i)), \circ)$ where $T2_{np^3,p}(C_{np^3}(R^{np^3,x+yp}_i))$ = $\{\theta_{np^3,p,jn}(C_{np^3}(R^{np^3,x+yp}_i))$ = $C_{np^3}(R^{np^3,x+yp}_{i+j}) :$ $j$ = $0,1,...,p-1$ and $i+j$ in $C_{np^3}(R^{np^3,x+yp}_{i+j})$ is calculated under addition modulo $p \}$, $1 \leq x \leq p-1$, $0 \leq y \leq np - 1$, $1 \leq x+yp \leq np^2-1$, $y\in\mathbb{N}_0$, $p,np^3-p\in R^{np^3,x+yp}_i$ and $i,n,x\in\mathbb{N}$. Using Theorem \ref{c4}, a list of $T2_{np^3,p}(C_{np^3}(R^{np^3,x+yp}_i))$, each containing $p$ isomorphic circulant graphs $C_{np^3}(R^{np^3,x+yp}_i)$ of Type-2 w.r.t. $m$ = $p$, for $p$ = 3,5,7, $n$ = 1,2 and $y$ = 0 is given in the Annexure and more such families of Type-2 isomorphic circulant graphs are presented in \cite{v24}.
\end{abstract}

\maketitle

	
\section{Introduction} This study is the $10^{th}$ part of a detailed study on Type-2 isomorphic circulant graphs containing ten parts by the authors \cite{v2-1}-\cite{v2-10}. Vilfred \cite{v96} defined Type-2 isomorphism, a new type of isomorphism different from already known circulant graph isomorphism of Adam's, of circulant graph $C_n(R)$ w.r.t. $m$ where $m$ = $\gcd(n, r) > 1$, $r\in R$ and $r,n\in\mathbb{N}$ and studied these graphs for $m$ = 2 in \cite{v13,v20}. Vilfred and Wilson \cite{v20}, \cite{vw1} - \cite{vw3} obtained Type-2 isomorphic circulant graphs $C_n(R)$ w.r.t. $m$ = 3,5,7 where $n\in\mathbb{N}$. We call Adam's isomorphism of circulant graphs as {\em Type-1 isomorphism} of circulant graphs. In \cite{v2-1}, Vilfred extended the definition of Type-2 isomorphism of circulant graphs $C_n(R)$ w.r.t. $m$ by considering $m > 1$ as a divisor of $\gcd(n, r)$ and $r\in R$. Here, $m$ = $\gcd(n, r) > 1$ and $r\in R$ is a particular case of $m > 1$ as a divisor of $\gcd(n, r)$ and $r\in R$. 

Based on the structures of families of isomorphic circulant graphs $C_n(R)$ of Type-2 w.r.t. $m$ = 2,3,5,7 obtained in \cite{v13},\cite{v20},\cite{vw1}-\cite{vw3}, we could find the structure of families of isomorphic circulant graphs $C_{np^3}(S)$ of Type-2 w.r.t. $m$ = $p$ and develope its theory that includes its related Abelian groups which we present in this paper where $p$ is a prime number and $n\in\mathbb{N}$. Theorems \ref{c1} and \ref{c4} are the main results and Theorems \ref{c10} and \ref{c12} are their more generalisation. In Theorem \ref{c4}, it is proved that for $i$ = 1 to $p$, circulant graphs $C_{np^3}(R^{np^3,x+yp}_i)$ are isomorphic of Type-2 w.r.t. $m$ = $p$ and they form Abelian group $(T2_{np^3,p}(C_{np^3}(R^{np^3,x+yp}_i)), \circ)$ where $T2_{np^3,p}(C_{np^3}(R^{np^3,x+yp}_i))$
= $\{\theta_{np^3,p,jn}(C_{np^3}(R^{np^3,x+yp}_i))$ = $C_{np^3}(R^{np^3,x+yp}_{i+j}) :$ $j$ = $0,1,...,p-1$ and $i+j$ in $C_{np^3}(R^{np^3,x+yp}_{i+j})$ is calculated under addition modulo $p \}$, $1 \leq x \leq p-1$, $0 \leq y \leq np - 1$, $1 \leq x+yp \leq np^2-1$, $y\in\mathbb{N}_0$, $p,np^3-p\in R^{np^3,x+yp}_i$ and $i,n,x\in\mathbb{N}$. Using Theorem \ref{c4}, a list of $T2_{np^3,p}(C_{np^3}(R^{np^3,x+yp}_i))$, each containing $p$ isomorphic circulant graphs $C_{np^3}(R^{np^3,x+yp}_i)$ of Type-2 w.r.t. $m$ = $p$, for $p$ = 3,5,7, $n$ = 1,2 and $y$ = 0 is given in the Annexure and more such families of Type-2 isomorphic circulant graphs are presented in \cite{v24}, $1 \leq i \leq p$.

This paper contains 4 sections and an annexure. Section 1 is introduction; Section 2 contains some basic definitions and results that are required in the subsequent sections; Section 3 presents families of isomorphic circulant graphs $C_{np^3}(R)$ of Type-2 w.r.t $m$ = $p$ $\ni$ $m > 1$ is a divisor of $\gcd(np^3, pq)$ = $p$, $r$ = $pq$ and $r\in C_{np^3}(R)$ and also present related Abelian groups where $p$ is a prime number, $n\in\mathbb{N}$ and $n \geq 2$ when $p$ = 2; Section 4 is the conclusion. It also contains an annexure containing a few families of Type-2 isomorphic circulant graphs. We follow remark \ref{r2.6} to establish Type-2 isomorphism w.r.t. $m$ among circulant graphs $C_n(R)$ and $C_n(S)$. Type-2 isomorphic circulant graphs don't have  Cayley Isomorphism (CI) property. For more on Type-2 circulant graph isomorphism one can refer \cite{v24, v2-1}. 

\section{Preliminaries }  

 Part 1 \cite{v2-1} contains many definitions and results related to Type-2 isomorphism of circulant graphs. In this section, we present a few definitions and results which are reqired in the subsequent sections.  
 
 \begin{definition}{\rm\cite{ad67}} \quad \label{d1} For $R =$ $\{r_1$, $r_2$, $\dots$, $r_k\}$ and $S$ = $\{s_1$, $s_2$, $\dots$, $s_k\}$, circulant graphs $C_n(R)$ and $C_n(S)$ are {\em Adam's isomorphic} or {\em Type-1 isomorphic} if there exists a positive integer $x$ $\ni$ $\gcd(n, x)$ = 1 and $S$ = $\{xr_1$, $xr_2$, $\dots$, $xr_k\}_n^*$ where $<r_i>_n^*$, the {\em reflexive modular reduction of a sequence} $< r_i >$, is the sequence obtained by reducing each $r_i$ under modulo $n$ to yield $r_i'$ and then replacing all resulting terms $r_i'$ which are larger than $\frac{n}{2}$ by $n-r_i'.$  
 \end{definition}
 
 \begin{definition}{\rm \cite{v2-1}} \label{d2} Let $Ad_n = \{\varphi_{n,x}: x\in \varphi_n\}$, $Ad_n(S) = \{\varphi_{n,x}(S): x\in \varphi_n\}$ = $\{xS: x\in \varphi_n\}$, $Ad_{n,x}(C_n(R))$ = $T1_{n,x}(C_n(R))$ = $\varphi_{n,x}(C_n(R))$ = $C_n(\varphi_{n,x}(R))$ = $C_n(xR)$, $x\in \varphi_n$ and $Ad_n(C_n(R)) = T1_n(C_n(R)) = \{\varphi_{n,y}(C_n(R)) = C_n(yR): y\in \varphi_n\}$ for sets $R,S \subseteq \mathbb{Z}_n$ where $\varphi_{n,x}(R)$ in $C_n(\varphi_{n,x}(R))$ is calculated under the reflexive modulo $n$. Define $'\circ'$ in $Ad_n(C_n(R))$ such that $\varphi_{n,x} \circ \varphi_{n,y}$ = $\varphi_{n,xy}$, $C_n(xR) \circ C_n(yR)$ = $C_n((xy)R)$ and $\varphi_{n,x}(C_n(R)) \circ \varphi_{n,y}(C_n(R))$ = $\varphi_{n,xy}(C_n(R))$, $\forall$ $x,y\in\varphi_n$.

Here, $(\varphi_{n,x} \circ \varphi_{n,y})(C_n(R))$ = $\varphi_{n,xy}(C_n(R))$ = $C_n((xy)R)$ = $C_n(xR) \circ C_n(yR)$  = $\varphi_{n,x}(C_n(R)) \circ \varphi_{n,y}(C_n(R))$, $\forall$ $x,y \in \varphi_n$. 
\end{definition}
  	
 Clearly, $Ad_n(C_n(R))$ is the set of all circulant graphs which are Adam's isomorphic to $C_n(R)$ and we call it as the {\em Adam's set} or {\em Type-1 set} of $C_n(R)$. Also, $(Ad_n(C_n(R)), \circ )$ = $(T1_n(C_n(R)), \circ )$ is an Abelian group and we call it as the {\em Adam's group} or {\em Type-1 group} of $C_n(R)$ under $'\circ'$. 
 
\begin{theorem} {\rm \cite{v2-1}}\quad  \label{t2.3} {\rm  Let $Ad_n(C_n(R))$ = $T1_n(C_n(R))$ = $\{\varphi_{n,x}(C_n(R)) = C_n(xR): x\in\varphi_n \}$. Then, $C_n(S)\in Ad_n(C_n(R))$ if and only if $Ad_n(C_n(R))$ = $Ad_n(C_n(S))$ if and only if $C_n(R)\in Ad_n(C_n(S))$. \hfill $\Box$  } 
 \end{theorem}

 In 1996, Vilfred \cite{v96} defined Type-2 isomorphism and studies on Type-2 isomorphism are found in \cite{v96}-\cite{v2-7}. In \cite{v2-1}, definition of Type-2 isomorphism w.r.t. $m$ of circulant graph $C_n(R)$ is extended by considering $m > 1$ as a divisor of $\gcd(n, r)$ and $r\in R$ instead of $m$ = $\gcd(n, r) > 1$ and $r\in R$.  Here, $m$ = $\gcd(n, r) > 1$ and $r\in R$ is a particular case of $m > 1$ is a divisor of $\gcd(n, r)$ and $r\in R$.

 \begin{definition} \cite{v2-1}  \quad  \label{d3} Let $V(K_n) = \{u_0,u_1,u_2,...,u_{n-1}\}$, $V(C_n(R)) = \{v_0,v_1,v_2,...,$ $v_{n-1}\},$ $r\in R$, $m > 1$ be a divisor of $\gcd(n, r)$ and $|R| \geq 3$.  Define one-to-one mapping $\theta_{n,m,t} :$ $V(C_n(R)) \rightarrow V(K_n)$ such that $\theta_{n,m,t}(v_x)$ = $u_{x+jtm}$,  $\theta_{n,m,t}((v_x, v_{x+s}))$ = $(\theta_{n,m,t}(v_x), \theta_{n,m,t}(v_{x+s}))$ under subscript arithmetic modulo $n$ and $\theta_{n,m,t}(C_n(R))$ = $C_n(\theta_{n,m,t}(R))$ for every $x \in \mathbb{Z}_n$, $s\in R$, $x$ = $qm+j$, $0 \leq j \leq m-1$, $0 \leq q,t \leq \frac{n}{m} -1$ and $\theta_{n,m,t}(R)$ in $C_n(\theta_{n,m,t}(R))$ is calculated under the reflexive modulo $n$. And for a particular value of $t,$ if  $\theta_{n,m,t}(C_n(R))$ = $C_n(S)$ for some $S$  and  $S \neq yR$ for all $y\in \varphi_n$ under reflexive modulo $n$, then $C_n(R)$ and $C_n(S)$ are called {\em isomorphic circulant graphs of Type-2 w.r.t. $m$} and the isomorphism as {\em Type-2 isomorphism w.r.t. $m$.} 

When $C_n(R)$ and $C_n(S)$ are Type-2 isomorphic w.r.t. $m$, then we also say that $C_{kn}(kR)$ $(= k.C_n(R))$ and $C_{kn}(kS)$ $(= k.C_n(S))$ are Type-2 isomorphic w.r.t. $m$. 	 
 \end{definition}
 
 \begin{rem}\quad \label{r2.5} Following steps are used to establish isomorphism of Type-2 w.r.t. $m$ between circulant graphs $C_n(R)$ and $C_n(S)$. (i) $R \neq S$ and $|R|$ = $|S| \geq 3$; (ii) $\exists$ $r\in R,S$ and $m > 1$ $\ni$ $m$ is a divisor of $\gcd(n, r)$ and for some $t$ $\ni$ $1 \leq t \leq \frac{n}{m} -1$, $\theta_{n,m,t}(C_n(R))$ = $C_n(S)$ and (iii) $S \neq xR$ for all $x\in\varphi_n$ under arithmetic reflexive modulo $n$. 
 \end{rem} 
 
 \begin{rem} \label{r2.6} \quad While searching for possible value(s) of $t$ for which the transformed graph $\theta_{n,m,t}(C_n(R))$ is circulant of the form $C_n(S)$ for some $S \subseteq [1, \frac{n}{2}],$ the calculation on $r_i$s which are integer multiples of $m$ need not be done  under the transformation $\theta_{n,m,t}$ as there is no change in these $r_i$s where $m > 1$ is a divisor of $\gcd(n, r)$ and $r\in R$. Also, for a given circulant graph $C_n(R)$, w.r.t. different values of $m$ or $r$ or both, we may get different Type-2 isomorphic circulant graphs.
 \end{rem}
 
Let $s \in \mathbb{Z}_n$,  $V_{n,m}$ = $\{\theta_{n,m,t}:$ $t = 0,1,...,\frac{n}{m}-1\}$, $V_{n,m}(s)$ = $\{\theta_{n,m,t}(s): t = 0,1,...,\frac{n}{m}-1\}$ and $V_{n,m}(C_n(R))$ = $\{\theta_{n,m,t}(C_n(R)): t = 0,1,...,\frac{n}{m}-1\}$. Define $'\circ'$ in $V_{n,m}$ such that $\theta_{n,m,t} ~\circ ~ \theta_{n,m,t'}$ =  $\theta_{n,m,t+t'}(\theta_{n,m,t} ~\circ ~ \theta_{n,m,t'})(x)$  $( = \theta_{n,m,t}(\theta_{n,m,t'}(x))$ = $\theta_{n,m,t}(x+jt'm)$ = $(x+jt'm)+jtm$ = $x+j(t+t')m )$ = $\theta_{n,m,t+t'}(x)$ and $\theta_{n,m,t}(C_n(R)) ~\circ ~ \theta_{n,m,t'}$ $(C_n(R))$ = $\theta_{n,m,t+t'}(C_n(R))$ for every $\theta_{n,m,t},\theta_{n,m,t'}\in V_{n,m}$ where $t+t'$ is calculated under addition modulo ~$\frac{n}{m}$. Clearly, $(V_{n,m}(s),~ \circ)$ and $(V_{n,m}(C_n(R)),~ \circ)$ are Abelian  groups, $\forall s \in \mathbb{Z}_n$. 

 \begin{definition} {\rm \cite{v2-3}} \quad Let $T2_{n,m}(C_n(R))$ = $\{C_n(R)\}$ $\cup$ $\{C_n(S):$ $C_n(S)$ is Type-2 isomorphic of $C_n(R)$ w.r.t. $m\}$ where $r\in R$ and $m > 1$ is a divisor of $\gcd(n, r)$. We call $T2_{n,m}(C_n(R))$ as {\em the Type-2 set of $C_n(R)$
w.r.t. $m$}. 
\end{definition}

Clearly, $T2_{n,m}(C_n(R)) \subseteq V_{n,m}(C_n(R))$ and $T1_n(C_n(R))$ $\cap$ $T2_{n,m}(C_n(R))$ = $\{C_n(R)\}$. 

 \begin{theorem}{\rm \cite{v2-3}} \quad {\rm  Let $C_n(R)$ $\cong$ $C_n(S)$, $R \neq S$, $|R| = |S| \geq 3$, $r\in R,S$ and $m > 1$ is a divisor of $\gcd(n, r)$. Then, $C_n(S)\in$ $T2_{n,m}(C_n(R))$ if and only if  $T2_{n,m}(C_n(R))$ = $T2_{n,m}(C_n(S))$ if and only if  $C_n(R)\in T2_{n,m}(C_n(S))$.  \hfill $\Box$}
 \end{theorem}
 
 \begin{theorem}{\rm \cite{v24}  \label{t2.8}  Under the above definition of $`\circ'$, $(T2_{n,m}(C_n(R)),  \circ)$ is a subgroup of $(V_{n,m}(C_n(R)), \circ)$ where $r\in R$ and $m > 1$ is a divisor of $\gcd(n, r)$. \hfill $\Box$}
 \end{theorem}
 
 \begin{definition}{\rm \cite{v24}}\quad \label{d2.9} With usual notation, group $(T2_{n,m}(C_n(R)), \circ)$ is called the {\em Type-2 group of $C_n(R)$} w.r.t.  $m$.  
 \end{definition}
 
\begin{definition} \cite{v24}\quad \label{d2.10} A circulant graph $C_n(R)$ is said to have {\em Cayley Isomorphism (CI) property} if whenever $C_n(S)$ is isomorphic to $C_n(R)$, they are Adam’s isomorphic. 
\end{definition}

CI-problem determines which graphs (or which groups) have the CI-property. Classification of cyclic CI-groups was completed by Muzychuk \cite{mu04} but investigation of circulant graphs without CI-property is not done much. Clearly, Type-2 isomorphic circulant graphs are circulant graphs without CI-property.  

\begin{definition}{\rm \cite{v13}}\quad \label{a14} For a set $R = \{r_1,r_2,\dots,r_k\}$ in $C_n(R)$, the {\it symmetric equidistance condition} w.r.t. $v_i$ is that  $v_{i}$ is adjacent to $v_{i+j}$  if and only if $v_{i}$ is adjacent to $v_{n-j+i}$, using subscript arithmetic modulo $n,$ $0 \leq i,j \leq n-1$. 
\end{definition} 

\begin{theorem}{\rm \cite{v13} \quad \label{ab14} For a set $R$ = $\{r_1,r_2,\dots,r_k\}  \subseteq [1, \frac{n}{2}]$, $r\in R$, $m > 1$ is a divisor of $\gcd(n, r)$, $1 \leq i \leq k$ and $0 \leq t \leq \frac{n}{m}-1$, $\theta_{n,m,t}(C_n(R))$ = $C_n(S)$ for some $S \subseteq [1, \frac{n}{2}]$ if and only if $\theta_{n,m,t}(C_n(R))$ satisfies the symmetric equidistance condition w.r.t. $v_0$.  \hfill $\Box$}
\end{theorem}

\begin{theorem}{\rm \cite{v20}}\quad \label{a17c} {\rm For $n \geq 2$, $1 \leq 2s-1 \leq 2n-1$, $n \neq 2s-1$, $R$ = $\{2,2s-1, 4n-(2s-1)\}$ and $S$ = $\{ 2,$ $2n-(2s-1),$ $2n+2s-1 \},$ $\theta_{8n,2,n}(C_{8n}(R))$ = $C_{8n}(S)$ = $\theta_{8n,2,3n}(C_{8n}(R))$, $\theta_{8n,2,n}(C_{8n}(S))$ = $C_{8n}(R)$ = $\theta_{8n,2,3n}(C_{8n}(S))$ and circulant graphs $C_{8n}(R)$ and $C_{8n}(S)$ are Type-2 isomorphic  w.r.t. $m$ = 2. When $n$ = $2s-1$, the two circulant graphs are the same.   \hfill $\Box$}
\end{theorem}

\begin{theorem} \label{a17d} {\rm \cite{v2-1} Let $n \geq 2$, $k \geq 3$, $1 \leq 2s-1 \leq 2n-1$, $n \neq 2s-1$, $R$ = $\{ 2s-1$, $4n-(2s-1)$, $2p_1$, $2p_2$, $\dots$, $2p_{k-2} \}$, $S$ = $\{2n-(2s-1)$, $2n+2s-1$, $2p_1,2p_2,\dots,2p_{k-2}\}$, $2y\in R,S$, $\gcd(4n,y)$ = 1, $p_1,p_2,\dots,p_{k-2} \in \mathbb{N}$ and $\gcd(p_1,p_2,\dots,p_{k-2})$ = 1. Then, (i) $\theta_{8n,2,n}(C_{8n}(R))$ = $C_{8n}(S)$ = $\theta_{8n,2,3n}(C_{8n}(R))$, $\theta_{8n,2,n}(C_{8n}(S))$ = $C_{8n}(R)$ = $\theta_{8n,2,3n}(C_{8n}(S))$ and (ii) corresponding to
each value of $k \geq 3$ and for a given set of values of $p_1,p_2,\dots,p_{k-2}$ and $n$, $C_{8n}(R)$ and $C_{8n}(S)$ are isomorphic of either Adam's or Type-2 w.r.t. $m$ = 2. Moreover, corresponding to each value of $k \geq 3$ and for all such possible values of $p_1,p_2,\dots,p_{k-2}$ and $n$, the set $\{ C_n(S)$ = $\theta_{8n,2,n}(C_{8n}(R)):$ $p_1,p_2,\dots,p_{k-2} \in \mathbb{N}\}$ contains all isomorphic circulant graphs of $C_{8n}(R)$ of Type-2 w.r.t. $m$ = 2.  \hfill $\Box$}
\end{theorem}

\begin{theorem}{\rm \cite{v20}}\quad \label{a18} {\rm For $n \geq 2$, $1 \leq 2s-1 < 2s'-1 \leq [\frac{n}{2}]$, $0 \leq t \leq [\frac{n}{2}]$, $R$ = $\{2,2s-1, 2s'-1\}$ and $n,s,s'\in \mathbb{N}$, if $\theta_{n,2,t}(C_n(R))$ and $C_n(R)$ are  isomorphic circulant graphs of Type-2 w.r.t. $m$ = 2 for some $t$, then $n \equiv 0~(mod ~ 8)$, $2s-1+2s'-1$ = $\frac{n}{2}$, $2s-1 \neq \frac{n}{8}$, $t$ = $\frac{n}{8}$ or $\frac{3n}{8}$, $1 \leq 2s-1 \leq \frac{n}{4}$ and $n \geq 16$. In particular, when $R$ = $\{2, 2s-1, 4n-(2s-1)\}$, $S$ = $\{2, 2n-(2s-1), 2n+2s-1\}$, $n\geq 2$ and $n,s\in \mathbb{N}$, $\theta_{8n,2,n}(C_{8n}(R))$ = $C_{8n}(S)$ = $\theta_{8n,2,3n}(C_{8n}(R))$ and $C_{8n}(R)$ and $C_{8n}(S)$ are Type-2 isomorphic w.r.t. $m$ = 2.    \hfill $\Box$}
\end{theorem}

\begin{theorem}{\rm \cite{v2-1} \quad \label{a20} For $n \geq 2$, $k \geq 3$, $1 \leq 2s-1 \leq 2n-1$, $n \neq 2s-1$, $R$ = $\{2, 2s-1, 4n-(2s-1)\}$ and $S$ = $\{2, 2n-(2s-1), 2n+2s-1\}$, $C_{8n}(R)$ and $C_{8n}(S)$ are Type-2 isomorphic w.r.t. $m$ = 2,  $T2_{8n,2}(C_{8n}(R))$ = $T2_{8n,2}(C_{8n}(S))$ = $\{C_{8n}(R), C_{8n}(S)\}$ and $(T2_{8n,2}(C_{8n}(R)), \circ)$ = $(T2_{8n,2}(C_{8n}(S)), \circ)$ is a Type-2 group, $n,s\in\mathbb{N}$. \hfill $\Box$ }
\end{theorem}
 
 \begin{theorem} \cite{vw1} \label{c41} {\rm For $R$ = $\{1, 3, 9n-1, 9n+1\}$, $S$ = $\{3, 3n+1, 6n-1, 12n+1\}$, $T$ = $\{3, 3n-1, 6n+1$, $12n-1\}$ and $n\in\mathbb{N}$, $\theta_{27n,3,n}(C_{27n}(R))$ = $C_{27n}(S)$, $\theta_{27n,3,n}(C_{27n}(S))$ = $C_{27n}(T)$, $\theta_{27n,3,n}(C_{27n}(T))$ = $C_{27n}(R)$ and $C_{27n}(R)$, $C_{27n}(S)$ and $C_{27n}(T)$ are Type-2 isomorphic circulant graphs w.r.t. $m$ = 3. \hfill $\Box$}
 \end{theorem}

 \begin{theorem} \label{c42} {\rm Let $k \geq 3$, $R$ = $\{1, 9n-1, 9n+1, 3p_1, 3p_2, \dots, 3p_{k-2}\}$, $S$ = $\{3n+1, 6n-1, 12n+1,$ $3p_1, 3p_2, \dots, 3p_{k-2}\}$, $T$ = $\{3n-1, 6n+1, 12n-1, 3p_1, 3p_2, \dots, 3p_{k-2}\}$, $\gcd(p_1,p_2,...,p_{k-2})$ = 1 and $k,n,p_1,p_2,\dots,p_{k-2}\in\mathbb{N}$. Then, $(i)$  $\theta_{27n,3,n}(C_{27n}(R))$ = $C_{27n}(S)$, $\theta_{27n,3,n}(C_{27n}(S))$ = $C_{27n}(T)$ and $\theta_{27n,3,n}(C_{27n}(T))$ = $C_{27n}(R)$ and $(ii)$ corresponding to each value of $k \geq 3$ and for a given set of values of $p_1,p_2,...,p_{k-2}$ and $n$, $C_{27n}(R)$, $C_{27n}(S)$ and $C_{27n}(T)$ are either all Type-2 isomorphic w.r.t. $m$ = 3 or all Adam's isomorphic.}
\end{theorem}
\begin{proof} The result follows from Theorem \ref{c41}, Remark \ref{r2.6} and definition of Type-2 isomorphism. 
\end{proof}

 \begin{theorem} \cite{vw2} \label{c43} {\rm For $R_i$ = $\{5, d_i, 25n-d_i, 25n+d_i, 50n-d_i, 50n+d_i\}$, $d_i$ = $5n(i-1)+1$, $i,j$ = 1 to 5 and $n\in\mathbb{N}$, $\theta_{125n,5,jn}(C_{125n}(R_i))$ = $C_{125n}(R_{i+j})$ and $C_{125n}(R_i)$ are Type-2 isomorphic circulant graphs w.r.t. $m$ = 5 where $i+j$ in $R_{i+j}$ is calculated under addition modulo 5.  		\hfill $\Box$}
\end{theorem}

\begin{theorem} \label{c44} {\rm Let $k \geq 3$, $d_i$ = $5n(i-1)+1$, $1 \leq i \leq 5$, $R_i$ = $\{d_i, 25n-d_i, 25n+d_i, 50n-d_i$, $50n+d_i, 5p_1, 5p_2, . . . , 5p_{k-2}\}$,  $k,n,p_1,p_2,...,p_{k-2}\in\mathbb{N}$ and $\gcd(p_1,p_2,...,p_{k-2}) = 1$. Then, corresponding to each value of $k \geq 3$ and for a given set of values of $p_1,p_2,...,p_{k-2}$ and $n$, circulant graphs $C_{125n}(R_i)$ are either all Type-2 isomorphic w.r.t. $m$ = 5 or all Adam's isomorphic, $1 \leq i \leq 5$.}
\end{theorem}
\begin{proof} The result follows from Theorem \ref{c43}, Remark \ref{r2.6} and definition of Type-2 isomorphism. 
\end{proof}

\begin{theorem} \cite{vw3} \label{c45} {\rm For $R_i$ = $\{7, d_i, 49n-d_i, 49n+d_i, 98n-d_i, 98n+d_i, 147n-d_i, 147n+d_i\}$, $d_i$ = $7n(i-1)+1$, $i,j$ = 1 to 7 and $n\in\mathbb{N}$, $\theta_{343n,7,jn}(C_{343n}(R_i))$ = $C_{343n}(R_{i+j})$ and $C_{343n}(R_i)$ are Type-2 isomorphic circulant graphs w.r.t. $m$ = 7 where $i+j$ is calculated under addition modulo 7. 	\hfill $\Box$}
\end{theorem}

\begin{theorem} \label{c46} {\rm Let $k \geq 3$, $d_i$ = $7n(i-1)+1$, $1 \leq i \leq 7$, $R_i$ = $\{d_i, 49n-d_i, 49n+d_i, 98n-d_i$, $98n+d_i$, $147n-d_i$, $147n+d_i$, $7p_1, 7p_2, . . . , 7p_{k-2}\}$, $k,n,p_1,p_2,...,p_{k-2}\in\mathbb{N}$ and $\gcd(p_1,p_2,...,p_{k-2}) = 1$. Then, corresponding to each value of $k \geq 3$ and for a given set of values of $p_1,p_2,...,p_{k-2}$ and $n$, circulant graphs $C_{343n}(R_i)$ are either all Type-2 isomorphic w.r.t. $m$ = 7 or all Adam's isomorphic, $1 \leq i \leq 7$.}
\end{theorem}
\begin{proof} The result follows from Theorem \ref{c45}, Remark \ref{r2.6} and definition of Type-2 isomorphism. 
\end{proof}

\section{Type-2 isomorphic circulant graphs $C_{np^3}(R)$ w.r.t. $m$ = $p$}  

Based on our studies on Type-2 isomorphic circulant graphs $C_n(R)$ w.r.t. $m$ = 2,3,5,7 in \cite{v20},\cite{vw1}-\cite{vw3}, we could find the structure of families of isomorphic circulant graphs $C_{np^3}(S)$ of Type-2 w.r.t. $m$ = $p$ and developed its theory that includes its related Abelian groups which are presented in this section where $p$ is a prime number and $n\in\mathbb{N}$ and $n \geq 2$ when $p$ = 2. Type-2 isomorphic circulant graphs are circulant graphs without the CI-property.  We use remark \ref{r2.5} to establish Type-2 isomorphism w.r.t.  $m$ among circulant graphs $C_n(R)$ and $C_n(S)$. 

\begin{theorem} \label{c1} {\rm Let $p$ be an odd prime number, $1 \leq i \leq p$, $1 \leq x \leq p-1$, $y\in\mathbb{N}_0$, $0 \leq y \leq np-1$, $1 \leq x+yp \leq np^2-1$, $d^{np^3, x+yp}_i = (i-1)xpn+x+yp$,  $R^{np^3, x+yp}_i$ $=$ $\{p$, $d^{np^3, x+yp}_i$, $np^2-d^{np^3, x+yp}_i$, $np^2+d^{np^3, x+yp}_i$, $2np^2-d^{np^3, x+yp}_i$, $2np^2+$ $d^{np^3, x+yp}_i,$ $3np^2-d^{np^3, x+yp}_i$, $3np^2+d^{np^3, x+yp}_i$, . . . , $(p-1)np^2$ - $d^{np^3, x+yp}_i$, $(p-1)np^2+d^{np^3, x+yp}_i$, $np^3-d^{np^3, x+yp}_i$, $np^3-p\}$ and $i,j,n,x\in\mathbb{N}$. Then, for a given set of values of $n$, $p$, $x$ and $y$, $\theta_{np^3,p,jn} (C_{np^3}(R^{np^3, x+yp}_i))$ = $C_{np^3}(R^{np^3, x+yp}_{i+j})$ and the $p$ circulant graphs $C_{np^3}(R^{np^3, x+yp}_i)$ are isomorphic of Type-2 w.r.t.  $p$, $1 \leq i,j \leq p$ where $i+j$ in $R^{np^3, x+yp}_{i+j}$ is calculated under addition modulo $p$ and $C_{np^3}(R^{np^3, x+yp}_0)$ = $C_{np^3}(R^{np^3, x+yp}_p)$.}
\end{theorem}
\begin{proof} We use remark \ref{r2.5} to establish Type-2 isomorphism w.r.t.  $m$ among circulant graphs $C_n(R)$ and $C_n(S)$. 

At first, let us prove $\theta_{np^3,p,jn} (R^{np^3, x+yp}_i)$ = $R^{np^3, x+yp}_{i+j}$ for $1 \leq i,j \leq p$, $1 \leq x \leq p-1$, $y\in\mathbb{N}_0$, $n\in\mathbb{N}$, $0 \leq y \leq np-1$ and $1 \leq x+yp \leq np^2-1$. 

For a given set of values of $n$, $x$ and $y$, we start with proving the above result for $i$ = 1 and 2. 

When $i$ = $1$, $n\in\mathbb{N}$, $y\in\mathbb{N}_0$, $1 \leq x \leq p-1$, $0 \leq y \leq np-1$ and $1 \leq x+yp \leq np^2-1$,

$d^{np^3, x+yp}_1 = x+yp$ and 

$R^{np^3, x+yp}_1$ = $\{p, x+yp, np^2-x-yp, np^2+x+yp, 2np^2-x-yp$, $2np^2+x+yp, 3np^2-x-yp$, 

\hfill $3np^2+x+yp$, $\dots$, $(p-1)np^2-x-yp$, $(p-1)np^2+x+yp$, $np^3-x-yp$, $np^3-p\}$.

When $i$ = $2$, $n\in\mathbb{N}$, $y\in\mathbb{N}_0$, $1 \leq x \leq p-1$, $0 \leq y \leq np-1$ and $1 \leq x+yp \leq np^2-1$, 

$d^{np^3, x+yp}_2 = xpn+x+yp$ and 

$R^{np^3, x+yp}_2 = \{p, xpn+x+yp, np^2-(xpn+x+yp), np^2+xpn+x+yp$,  

\hspace{2.5cm} $2np^2-(xpn+x+yp)$, $2np^2+xpn+x+yp$,   

\hspace{2.5cm} $3np^2-(xpn+x+yp)$, $3np^2+xpn+x+yp$, $\dots$,  

\hfill $(p-1)np^2-(xpn+x+yp)$, $(p-1)np^2+xpn+x+yp$, $np^3-(xpn+x+yp)$, $np^3-p\}$. 

For $1 \leq i,j \leq p$, $n\in\mathbb{N}$, $y\in\mathbb{N}_0$, $1 \leq x \leq p-1$, $0 \leq y \leq np-1$ and $1 \leq x+yp \leq np^2-1$, using the definition of $\theta_{n,m,t}$, we get, 
\\
$\theta_{np^3,p,n}(R^{np^3, x+yp}_1)$ = $\theta_{np^3,p,n}(\{p, np^3-p \})$  

\hfill $\bigcup \theta_{np^3,p,n}(\{x+yp, np^2+x+yp, 2np^2+x+yp, 3np^2+x+yp, \dots, (p-1)np^2+x+yp \})$ 

\hfill $ \bigcup \theta_{np^3,p,n}(\{np^2-x-yp, 2np^2-x-yp, 3np^2-x-yp, \dots, (p-1)np^2-x-yp,  np^3-x-yp \})$

= $\{p, np^3-p \}$ 

\hfill $ \bigcup (xpn+ \{x+yp, np^2+x+yp, 2np^2+x+yp, 3np^2+x+yp, \dots, (p-1)np^2+x+yp \})$ 

\hfill $\bigcup ((p-x)pn+ \{np^2-x-yp, 2np^2-x-yp, 3np^2-x-yp,$ 

\hfill $\dots, (p-1)np^2-x-yp, np^3-x-yp \})$ 

 = $\{p, np^3-p, xpn+x+yp, np^2+xpn+x+yp, 2np^2+xpn+x+yp, 3np^2+xpn+x+yp, \dots,$  

\hfill $  (p-1)np^2+xpn+x+yp, 2np^2-(xpn+x+yp), 3np^2-(xpn+x+yp),  \dots,$ 

\hfill $  np^3-(xpn+x+yp), np^2-(xpn+x+yp)\} = R^{np^3,x+yp}_2;$  

 $\theta_{np^3,p,in}(R^{np^3, x+yp}_1)$ = $\theta_{np^3,p,in}(\{p, np^3-p \})$ $\bigcup \theta_{np^3,p,in}(\{x+yp, np^2+x+yp$, $2np^2+x+yp$, 

\hfill $3np^2+x+yp, \dots, (p-1)np^2+x+yp \}) \bigcup \theta_{np^3,p,in}(\{np^2-x-yp$,  

\hfill $2np^2-x-yp$, $3np^2-x-yp, \dots, (p-1)np^2-x-yp,  np^3-x-yp \})$

= $\{p, np^3-p \} \bigcup (ixpn+ \{x+yp, np^2+x+yp, 2np^2+x+yp, 3np^2+x+yp, \dots$, 

\hfill $ (p-1)np^2+x+yp \}) \bigcup ~(i(p-x)pn+ \{np^2-x-yp, 2np^2-x-yp, 3np^2-x-yp, \dots$, 

~\hfill $(p-1)np^2-x-yp, np^3-x-yp \})$ 

 = $\{p, np^3-p, ixpn+x+yp, np^2+ixpn+x+yp, 2np^2+ixpn+x+yp, 3np^2+ixpn+x+yp, \dots,$ 

~ \hfill $(p-1)np^2+ixpn+x+yp, (1+i)np^2-(ixpn+x+yp),$ $(2+i)np^2-(ixpn+x+yp),$  

~\hfill  $(3+i)np^2-(ixpn+x+yp), \dots,$ $(p-1+i)np^2-(ixpn+x+yp),$

~\hfill $(p+i)np^2-(ixpn+x+yp)\} = R^{np^3,x+yp}_{i+1}$  \\
since $\{1+i, 2+i, 3+i, . . . , (p-1)+i, p+i\}$ = $\{1, 2, . . . , p-1, p = 0\}$ under addition modulo $p$ and $d^{np^3,x+yp}_{i+1} = ixpn+x+yp$, $1 \leq i \leq p$. Thus, the above result is true for $i$ = 1.

In a similar way, we can prove the above result for $i$ = 2 and also for the general case that $\theta_{np^3,p,jn}(R^{np^3,x+yp}_i)$ = $R^{np^3,x+yp}_{i+j}$ where $i+j$ in $R^{np^3,x+yp}_{i+j}$ is calculated under addition modulo $p$, $1 \leq x \leq p-1$, $0 \leq y \leq np-1$, $1 \leq x+yp \leq np^2-1$, $y\in\mathbb{N}_0$, $n\in\mathbb{N}$ and $1 \leq i,j \leq p$. And thereby, we get, $\theta_{np^3,p,jn}(C_{np^3}(R^{np^3,x+yp}_i))$ = $C_{np^3}(R^{np^3,x+yp}_{i+j})$ where $i+j$ in $R^{np^3,x+yp}_{i+j}$ is calculated under addition modulo $p$, $1 \leq x \leq p-1$, $0 \leq y \leq np-1$, $1 \leq x+yp \leq np^2-1$, $y\in\mathbb{N}_0$, $n\in\mathbb{N}$ and $1 \leq i,j \leq p$. 

From the definition of $\theta_{n,m,t}$ acting on $C_{n}(R)$, we get, for a given set of values of $n$, $x$ and $y$, the circulant graphs $C_{np^3}(R^{np^3, x+yp}_i)$ are isomorphic, $1 \leq x \leq p-1$, $0 \leq y \leq np-1$, $1 \leq x+yp \leq np^2-1$, $y\in\mathbb{N}_0$, $n\in\mathbb{N}$ and $1 \leq i,j \leq p$. To complete the proof, we have to establish their Type-2 isomorphism. 

\vspace{.2cm}
\noindent
{\bf {\it Claim.}} For a given set of values of $n$, $x$ and $y$, $C_{np^3}(R^{np^3, x+yp}_i)$ are Type-2 isomorphic w.r.t.  $p$, $1 \leq i \leq p$, $n\in\mathbb{N}$, $y\in\mathbb{N}_0$, $1 \leq x \leq p-1$, $0 \leq y \leq np-1$ and $1 \leq x+yp \leq np^2-1$.  

 For $1 \leq i,j\leq p$, $1 \leq x \leq p-1$, $1 \leq x+yp \leq np^2-1$, $0 \leq y \leq np-1$, $y\in\mathbb{N}_0$ and $n\in\mathbb{N}$ and given $n$, $x$ and $y$, $d^{np^3,x+yp}_i = (i-1)xpn+x+yp\in$ $R^{np^3,x+yp}_i$. And $d^{np^3,x+yp}_i$ = $d^{np^3,x+yp}_j$ if and only if $(i-1)xpn+x+yp$ = $(j-1)xpn+x+yp$ if and only if $i$  = $j$ if and only if $R^{np^3, x+yp}_i$ = $R^{np^3, x+yp}_j$. Thus, for given $n$, $x$ and $y$ and different $i,$ all the $p$ sets $R^{np^3, x+yp}_i$ are different and thereby all the $p$ circulant graphs $C_{np^3}(R^{np^3, x+yp}_i)$ are also distinct, $1 \leq i,j\leq p$, $1 \leq x \leq p-1$, $0 \leq y \leq np-1$, $y\in\mathbb{N}_0$, $n\in\mathbb{N}$ and $1 \leq x+yp \leq np^2-1$. We have already proved that for given $n$, $x$ and $y$, $\theta_{np^3,p,jn}(C_{np^3}(R^{np^3, x+yp}_i))$ = $C_{np^3}(R^{np^3, x+yp}_{i+j})$ and thereby all the $p$ circulant graphs $C_{np^3}(R^{np^3, x+yp}_i)$ are isomorphic where $i+j$ in $R^{np^3, x+yp}_{i+j}$ is calculated under addition modulo $p$. This implies that all the $p$ circulant graphs $C_{np^3}(R^{np^3, x+yp}_i)$ are distinct but isomorphic, $1 \leq i\leq p$, $1 \leq x \leq p-1$, $0 \leq y \leq np-1$, $1 \leq x+yp \leq np^2-1$, $y\in\mathbb{N}_0$ and $n\in\mathbb{N}$.

To prove their Type-2 isomorphism w.r.t.  $p$, it is enough to prove, each pair of circulant graphs $C_{np^3}(R^{np^3, x+yp}_i)$ and $C_{np^3}(R^{np^3, x+yp}_j)$ for $i \neq j$ are not of Type-1 for a fixed $n$, $x$ and $y$ where $1 \leq i,j\leq p$, $1 \leq x \leq p-1$, $0 \leq y \leq np-1$, $1 \leq x+yp \leq np^2-1$, $y\in\mathbb{N}_0$ and $n\in\mathbb{N}$. We start with the circulant graph $C_{np^3}(R^{np^3, x+yp}_1)$. 

\vspace{.2cm}
\noindent
{\bf{\it Sub-claim.}} $C_{np^3}(R^{np^3, x+yp}_1)$ and $C_{np^3}(R^{np^3, x+yp}_i)$ are Type-2 isomorphic w.r.t.  $p$ for given $n$, $x$ and $y$ and $2 \leq i \leq p$.

\vspace{.2cm}
If not, they are of Adam's isomorphic. This implies, there exists $s,q \in \mathbb{N}$ such that $C_{np^3}(sR^{np^3, x+yp}_1)$ = $C_{np^3}(R^{np^3, x+yp}_i)$ where $s$ = $qpn-j,$ $\gcd(np^3, s)$ = $1$, $1 \leq j,x \leq p-1,$ $1 \leq qpn-j \leq np^3-1$, $n\in\mathbb{N}$, $2 \leq i \leq p$, $0 \leq y \leq np-1$, $1 \leq x+yp \leq np^2-1$ and $y\in\mathbb{N}_0$. This also implies, $1 \leq j \leq p-1$ and $\gcd(qn, j)$ = 1. Now, consider the case when $j$ = $1.$ In this case, $s$ = $qpn-1,$ $\gcd(np^3, qpn-1)$ = $1,$ $\gcd(qn, 1)$ = 1, $1 \leq qpn-1 \leq np^3-1,$ $C_{np^3}((qpn-1)R^{np^3,x+yp}_1)$ = $C_{np^3}(R^{np^3,x+yp}_i),$ $2 \leq i \leq p$, $1 \leq x \leq p-1$, $0 \leq y \leq np-1$,  $1 \leq x+yp \leq np^2-1$, $y\in\mathbb{N}_0$, $n\in\mathbb{N}$. This implies, $(qpn-1)\{p,$ $x+yp,$ $np^2-x-yp$, $np^2+x+yp,$ $2np^2-x-yp,$ $2np^2+x+yp,$ $3np^2-x-yp,$ $3np^2+x+yp,$ $\dots,$ $(p-1)np^2-x-yp,$ $(p-1)np^2+x+yp,$ $np^3-x-yp,$ $np^3-p\}$ = $\{p,$ $(i-1)xpn+x+yp,$ $np^2-((i-1)xpn+x+yp),$ $np^2+(i-1)xpn+x+yp,$ $2np^2-((i-1)xpn+x+yp),$ $2np^2+(i-1)xpn+x+yp,$ $3np^2-((i-1)xpn+x+yp),$ $3np^2+(i-1)xpn+x+yp,$ $\dots,$ $(p-1)np^2-((i-1)xpn+x+yp),$ $(p-1)np^2+(i-1)xpn+x+yp,$ $np^3-((i-1)xpn+x+yp),$ $np^3-p\}$ under arithmetic modulo $np^3,$ $2 \leq i \leq p$, $1 \leq x \leq p-1$, $y\in\mathbb{N}_0$, $n\in\mathbb{N}$, $0 \leq y \leq np-1$, $1 \leq x+yp \leq np^2-1$. This implies, $(qpn-1)p,$ $(qpn-1)(np^3-p),$ $p+p_1np^3$ and $np^3-p+p_2np^3$ are the only numbers, each is a multiple of $p,$ in the two sets for some $p_1,p_2\in \mathbb{N}_0,$ $2 \leq i \leq p$, $1 \leq x \leq p-1$, $1 \leq qpn-1 \leq np^3-1$, $\gcd(np^3, qpn-1)$ = $1$, $0 \leq y \leq np-1$, $1 \leq x+yp \leq np^2-1$, $y\in \mathbb{N}_0$ and $n,q\in\mathbb{N}$. Under this, the following two cases arise.

\vspace{.2cm}
\noindent
{\bf {\it Case 1.}} $(qpn-1)p = p+p_1np^3$, $p_1 \in \mathbb{N}_0,$~ $n,q \in \mathbb{N},$ $1 \leq qpn-1 \leq np^3-1.$ 

In this case, the possible values of $p_1$ are $~0,~1,~2,~\dots,~p-1~$ since $1 \leq qpn-1 \leq np^3-1$ and $q \in \mathbb{N}.$ When $p_1$ = $0$, $qpn-1$ = $1$ and $\gcd(qn, j)$ = 1; $p_1 = 1$, $qpn-1$ = $np^2+1$ and $\gcd(qn, j)$ = 1; $p_1$ = $2$, $qpn-1$ = $2np^2+1$ and $\gcd(qn, j)$ = 1; $\dots;$ $p_1$ = $p-1$, $qpn-1$ = $(p-1)np^2+1$ and $\gcd(qn, j)$ = 1. Now, let us calculate $(qpn-1)R^{np^3,x+yp}_1$ for $qpn-1$ = $np^2+1,$ $2np^2+1,$ $\dots,$ $(p-1)np^2+1$ under arithmetic modulo $np^3.$ 

When $qpn-1$ = $np^2+1,$ for given $n$, $x$ and $y$, $1 \leq x \leq p-1$, $y\in\mathbb{N}_0$, $0 \leq y \leq np-1$, $1 \leq x+yp \leq np^2-1$ and $n,x\in\mathbb{N}$, under arithmetic modulo $np^3,$ 

\vspace{.2cm}
\noindent
$(qpn-1)R^{np^3,x+yp}_1$ = $(np^2+1)R^{np^3,x+yp}_1$ 

 = $(np^2+1)\{p, x+yp, np^2-x-yp, np^2+x+yp, 2np^2-x-yp, 2np^2+x+yp, $  

 \hfill $3np^2-x-yp, 3np^2+x+yp, \dots, $ $(p-1)np^2-x-yp, (p-1)np^2+x+yp, np^3-x-yp, np^3-p\}$ 

 = $\{p, xnp^2+x+yp, (p-x+1)np^2-x-yp, (x+1)np^2+x+yp, (p-x+2)np^2-x-yp$, 

\hfill   $(p-x+2)np^2-x-yp, (x+2)np^2+x+yp, (p-x+3)np^2-x-yp, (x+3)np^2+x+yp, \dots,$ 

~\hfill  $(p-x+p-1)np^2-x-yp, (x-1)np^2+x+yp, (p-x)np^2-x-yp, np^3-p\} = R^{np^3,x+yp}_1.$ 

\vspace{.2cm}
\noindent
Here, $\{xnp^2+x+yp$, $(x+1)np^2+x+yp$, $(x+2)np^2+x+yp$, $(x+3)np^2+x+yp$, $\dots$, 

\hfill $(x-1)np^2+x+yp$ = $(x+(p-1))np^2+x+yp\}$ 

=  $\{x+yp$, $np^2+x+yp$, $2np^2+x+yp$, $3np^2+x+yp$, $\dots$, $(p-1)np^2+x+yp\}$ and 
\\
$\{(p-x+1)np^2-x-yp$, $(p-x+2)np^2-x-yp$, $(p-x+3)np^2-x-yp$, . . . , 

\hfill $(p-x+p-1)np^2-x-yp)$, $(p-x+p)np^2-x-yp)\}$ 

= $\{np^2-x-yp$, $2np^2-x-yp$, $3np^2-x-yp$, $\dots$, $(p-1)np^2-x-yp$, $np^3-x-yp\}$ under arithmetic modulo $np^3$, $1 \leq x \leq p-1$, $y\in\mathbb{N}_0$, $1 \leq x+yp \leq np^2-1$, $0 \leq y \leq np-1$ and $n,x\in\mathbb{N}$.

Similarly, we can prove that $(qpn-1)R^{np^3,x+yp}_1$ = $R^{np^3,x+yp}_1$ when $qpn-1$ = $2np^2+1,$ $3np^2+1,$ $\dots,$ $(p-1)np^2+1$ under arithmetic modulo $np^3.$ This implies, $C_{np^3}((qpn-1)R^{np^3,x+yp}_1)$ = $C_{np^3}(R^{np^3,x+yp}_1)$ and $\neq$ $C_{np^3}(R^{np^3,x+yp}_i)$ for $2 \leq i \leq p-1$, $qpn-1$ = $np^2+1,$ $2np^2+1,$ $\dots,$ $(p-1)np^2+1$, $1 \leq x \leq p-1$, $y\in\mathbb{N}_0$, $1 \leq x+yp \leq np^2-1$, $0 \leq y \leq np-1$ and $n,x\in\mathbb{N}$. 

In a similar way, we can prove that for $2 \leq j \leq p-1$ and $\gcd(qpn-j, np^3)$ = 1, $(qpn-j)$ $R^{np^3,x+yp}_1$ = $R^{np^3,x+yp}_1$ and $\neq$ $R^{np^3,x+yp}_i$ for $2 \leq i \leq p-1$, $y\in\mathbb{N}_0$, $1 \leq x \leq p-1$, $0 \leq y \leq np-1$, $1 \leq x+yp \leq np^2-1$, $n,x\in\mathbb{N}$ and when $nqp-j$ = $np^2+1,$ $2np^2+1,$ $\dots,$ $(p-1)np^2+1$. This implies, for $2 \leq i \leq p-1$ and given $n$, $x$ and $y$, $C_{np^3}(R^{np^3, x+yp}_1)$ and $C_{np^3}(R^{np^3, x+yp}_i)$ are not Adam's isomorphic, $1 \leq x \leq p-1$, $1 \leq x+yp \leq np^2-1$, $0 \leq y \leq np-1$, $y\in\mathbb{N}_0$ and $n,x\in\mathbb{N}$. Hence the sub-claim is true in this case. 

\vspace{.2cm}
\noindent
{\bf {\it Case 2.}} $(qpn-1)p = np^3-p+p_2np^3$, $p_2 \in \mathbb{N}_0$, $1 \leq qpn-1 \leq np^3-1$ and $n,q \in \mathbb{N}$. 

In this case, the possible values of $p_2$ are $0,~1,~2,~.~.~.~,~p-1$ since ~ $1 \leq qpn-1 \leq np^3-1$ and $n,q \in \mathbb{N}$. When $p_2$ = $0$, $qpn-1$ = $np^2-1$; $p_2$ = $1$, $qpn-1$ = $2np^2-1$; $\dots$; $p_2$ = $p-1$, $qpn-1$ = $np^3-1$. Now, let us calculate $(qpn-1)R^{np^3,x+yp}_1$ for $qpn-1$ = $np^2-1$, $2np^2-1$, $\dots$, $np^3-1$ under arithmetic modulo $np^3$. 

When $nqp-1$ = $np^2-1$, for given $n$, $x$ and $y$, $1 \leq x \leq p-1$, $y\in\mathbb{N}_0$, $0 \leq y \leq np-1$, $1 \leq x+yp \leq np^2-1$ and $n,x\in\mathbb{N}$, under arithmetic modulo $np^3$, 

\vspace{.2cm}
\noindent
$(qpn-1)R^{np^3,x+yp}_1 = (np^2-1)R^{np^3,x+yp}_1 $ 

 = $(np^2-1)\{p, x+yp, np^2-x-yp, np^2+x+yp, 2np^2-x-yp,2np^2+x+yp,  3np^2-x-yp$,  

 \hfill $3np^2+x+yp, \dots, (p-1)np^2-x-yp, (p-1)np^2+x+yp, np^3-x-yp, np^3-p\}$ 

   = $\{np^3-p, xnp^2-x-yp, (p-(x+1))np^2+x+yp$, $(x-1)np^2-x-yp, (p-(x+2))np^2+x+yp$, 

\hspace{.5 cm} $(x-2)np^2-x-yp, (p-(x+3))np^2+x+yp$, $(x-3)np^2-x-yp, \dots, (p-(x-1))np^2+x+yp$, 

~\hfill  $(p+x-(p-1))np^2-x-yp, (p-x)np^2+x+yp, p\}$ = $R^{np^3,x+yp}_1$. 

\vspace{.2cm}
Here, $\{xnp^2-x-yp$, $(x-1)np^2-x-yp$, $(x-2)np^2-x-yp$, $(x-3)np^2-x-yp$, ..., 

\hfill $(p+x-(p-1))np^2$ $-x-yp\}$ 

\hfill = $\{np^2-x-yp$, $2np^2-x-yp$, $3np^2-x-yp$, $\dots$, $(p-1)np^2-x-yp$, $np^3-x-yp\}$ and 
\\
$\{(p-(x+1))np^2+x+yp$, $(p-(x+2))np^2+x+yp$, $(p-(x+3))np^2+x+yp$, $\dots$, 

\hfill $(p-(x-1))np^2+x+yp)$, $(p-x)np^2+x+yp\}$ 

\hfill = $\{x+yp$, $np^2+x+yp$, $2np^2+x=yp$, $3np^2+x+yp$, $\dots$, $(p-1)np^2+x+yp\}$ under 
\\
arithmetic modulo $np^3$, $1 \leq x \leq p-1$, $y\in\mathbb{N}_0$, $0 \leq y \leq np-1$, $1 \leq x+yp \leq np^2-1$ and $n,x\in\mathbb{N}$.

Similarly, we can prove that $(qpn-1)R^{np^3,x+yp}_1$ = $R^{np^3,x+yp}_1$ when $qpn-1$ = $2np^2-1,$ $3np^2-1,$ $\dots,$ $np^3-1,$ under arithmetic modulo $np^3.$ This implies that $C_{np^3}((qpn-1)R^{np^3,x+yp}_1)$ = $C_{np^3}(R^{np^3,x+yp}_1)$ when $qpn-1$ = $np^2-1,$ $2np^2-1,$ $\dots,$ $np^3-1.$ In a similar way, we can prove that $(qpn-j)R^{np^3,x+yp}_1$ = $R^{np^3,x+yp}_1,$ under arithmetic modulo $np^3$, when $qpn-j$ = $np^2-1,$ $2np^2-1,$ $\dots,$ $np^3-1$ for $2 \leq j \leq p-1$ and $\gcd(qpn-j, np^3)$ = 1. This implies, for every $s = qpn-j$ which is relative prime to $np^3$ and for given $n$, $x$ and $y$, $C_{np^3}(sR^{np^3,x+yp}_1)$ = $C_{np^3}(R^{np^3,x+yp}_1)$ and $\neq$ $C_{np^3}(R^{np^3,x+yp}_i)$ for $2 \leq i \leq p-1$, $1 \leq j,x \leq p-1$,  $y\in\mathbb{N}_0$, $0 \leq y \leq np-1$, $1 \leq x+yp \leq np^2-1$ and for $nqp-j$ = $np^2-1,$ $2np^2-1,$ $\dots,$ $(p-1)np^2-1$, $1 \leq qpn-j \leq np^3-1$ and $n,q,x\in\mathbb{N}$. Hence the sub-claim is also true in this case. 

This implies, for given $n$, $x$ and $y$ and for $i = 2,3,\dots,p-1$, $C_{np^3}(R^{np^3,x+yp}_1)$ and $C_{np^3}(R^{np^3,x+yp}_i)$ can not be Type-1 isomorphic and so they are Type-2 isomorphic w.r.t.  $m$ = $p$ since $\theta_{np^3,p,jn} (C_{np^3}(R^{np^3, x+yp}_i))$ = $C_{np^3}(R^{np^3, x+yp}_{i+j})$ for $1 \leq i,j \leq p$, $1 \leq x \leq p-1$,  $y\in\mathbb{N}_0$, $n,x\in\mathbb{N}$, $0 \leq y \leq np-1$, $1 \leq x+yp \leq np^2-1$ and $i+j$ in $R^{np^3, x+yp}_{i+j}$ is calculated under arithmetic modulo p. Hence the claim is true.

Similarly, we can show that for given $n$, $x$ and $y$, $1 \leq i,j \leq p$ and $i \neq$ $j$, $C_{np^3}(R^{np^3,x+yp}_i)$ and $C_{np^3}(R^{np^3,x+yp}_j)$ can not be Adam's isomorphic by proving $C_{np^3}((qpn-t)R^{np^3,x+yp}_i)$ = $C_{np^3}(R^{np^3,x+yp}_i)$ when $qpn-t$ = $np^2-1$, $2np^2-1$, $\dots$, $np^3-1$ as well as $qpn-t$ = $np^2+1$, $2np^2+1$, $\dots$, $(p-1)np^2+1$ where $\gcd(np^3, qpn-t)$ = 1, $1 \leq qpn-t \leq np^3-1$, $1 \leq t,x \leq p-1$, $0 \leq y \leq np-1$, $1 \leq x+yp \leq np^2-1$, $y\in\mathbb{N}_0$ and $n,q,t,x\in\mathbb{N}$.   

This implies, for given $n$, $x$ and $y$, all the $p$ isomorphic circulant graphs $C_{np^3}(R^{np^3,x+yp}_i)$ are Type-2 isomorphic w.r.t.  $m$ = $p$ for $i$ = 1 to $p$, $y\in\mathbb{N}_0$, $1 \leq x \leq p-1$, $0 \leq y \leq np-1$, $1 \leq x+yp \leq np^2-1$ and $n,x\in\mathbb{N}$. Hence, we get the result.  
\end{proof}

\begin{theorem} \label{c2} {\rm Let $p$ be an odd prime number, $k \geq 3$, $1 \leq i,j \leq p$, $1 \leq x \leq p-1$, $y\in\mathbb{N}_0$, $0 \leq y \leq np - 1$, $1 \leq x+yp \leq np^2-1$, $d^{np^3, x+yp}_i$ = $(i-1)xpn+$ $x+yp$, $R^{np^3, x+yp}_i$ = $\{d^{np^3, x+yp}_i,$ $np^2-d^{np^3, x+yp}_i,$ $np^2+d^{np^3, x+yp}_i,$ $2np^2-$ $d^{np^3, x+yp}_i,$ $2np^2+d^{np^3, x+yp}_i,$ $3np^2-d^{np^3, x+yp}_i,$ $3np^2+d^{np^3, x+yp}_i,$ ..., $(p-1)np^2-d^{np^3, x+yp}_i,$ $(p-1)np^2+d^{np^3, x+yp}_i,$ $np^3-d^{np^3, x+yp}_i,$ $pp_1,$ $pp_2,$ $\dots$, $pp_{k-2},$ $p(np^3-p_{k-2})$, $p(np^3-p_{k-3})$, . . . , $p(np^3-p_1)\}$, $\gcd(p_1,p_2,...,p_{k-2}) = 1$ and $i,j,k,n,x,p_1,p_2,...,p_{k-2} \in \mathbb{N}.$  Then, corresponding to each value of $k \geq 3$ and for a given set of values of $p, x, y, p_1, p_2, ... , p_{k-3},p_{k-2}, n$ and for $i$ = $1,2,...,p$, $(i)$ $\theta_{np^3,p,jn}(C_{np^3}(R^{np^3, x+yp}_i))$ = $C_{np^3}(R^{np^3, x+yp}_{i+j})$ and $(ii)$ the $p$ circulant graphs $C_{np^3}(R^{np^3, x+yp}_i)$ are either all Type-1 or all Type-2 (and without $CI$-property) w.r.t.  $m$ = $p$ where $i+j$ in $R^{np^3, x+yp}_{i+j}$ is calculated under addition modulo $p$, $1 \leq j \leq p$.}
\end{theorem}
\begin{proof}\quad Using Theorem \ref{c1}, for given $n$, $x$ and $y$, the $p$ circulant graphs $C_{np^3}(R^{np^3, x+yp}_i)$ for $i$ = $1$ to $p$ are Type-2 isomorphic w.r.t.  $m$ = $p$ where $p$ is an odd prime, $d^{np^3, x+yp}_i = (i-1)xpn+x+yp$ and $R^{np^3, x+yp}_i$ = $\{p$, $d^{np^3, x+yp}_i$, $np^2-d^{np^3, x+yp}_i$, $np^2+d^{np^3, x+yp}_i$, $2np^2-d^{np^3, x+yp}_i$, $2np^2+$ $d^{np^3, x+yp}_i,$ $3np^2-d^{np^3, x+yp}_i$, $3np^2+d^{np^3, x+yp}_i,$ . . . , $(p-1)np^2-d^{np^3, x+yp}_i$, $(p-1)np^2+d^{np^3, x+yp}_i,$ $np^3-d^{np^3, x+yp}_i,$ $np^3-p\}$. Then the result follows from remark \ref{r2.6}, definition of Type-2 isomorphic circulant graphs and the property that Type-2 isomorphic circulant graphs are circulant graphs without $CI$-property. 
\end{proof}

Our next theorem is on $T2_{np^3,p}(C_{np^3}(R))$ and before that, let us consider a related problem. 

\begin{prm} \label{c3} {\rm For $R_1$ = $\{3,7,20,34,47,61,74,78\}$, $R_2$ = $\{3,11,16,38,43,65,70,78\}$ and 
		\\
		$R_3$ = $\{2,3,25,29,52,56,78,79\}$, prove the following.
		\begin{enumerate}
			\item [\rm (i)] $C_{81}(R_i)$ are Type-2 isomorphic w.r.t. $m$ = 3 for $i$ = 1,2,3;
			\item [\rm (ii)] Find $T2_{81,3}(C_{81}(R_i))$ for $i$ = 1,2,3; and
			\item [\rm (iii)] Find $(T2_{81,3}(C_{81}(R_i)), \circ)$ for $i$ = 1,2,3.
	\end{enumerate}	    } 
\end{prm}
\noindent
{\bf Solution.} All the three sets $R_1$, $R_2$ and $R_3$ are symmetric subsets of $\mathbb{Z}_{81}$, each contains 3 and 78 as common elements and $\gcd(81, 3)$ = 3 = $\gcd(81, 78)$ and hence $C_{81}(R_i)$ may be of Type-2 isomorphic w.r.t.  $m$ = 3, $i$ = 1,2,3. Let $m$ = 3 = $p$ = $n$ as in Theorem \ref{c1}. This implies, $np$ = 9, $np^2$ = 27, $np^3$ = 81 and $\gcd(np^3, 3)$ = 3 = $m$. 

The minimum jump size which is not a multiple of 3 in $R_1$ is 7 which implies, it is possible to consider 7 = $x+yp$ as in Theorem \ref{c1}. In this case, $x = 1$, $y = 2$ and so $d^{81,7}_1$ = $x+yp$ = 7, $d^{81,7}_2$ = $npx+x+yp$ = 9+7 = 16 and $d^{81,7}_3$ = $2npx+x+yp$ = 18+7 = 25. Using Theorem \ref{c1}, we get,

\vspace{.1cm}
$R^{81,7}_1$ = $\{3,7,20,34,47,61,74,78\}$ = $R_1$ = $\theta_{81,3,0}(R^{81,7}_1)$,  

\vspace{.1cm}
$R^{81,7}_2$ = $\{3,11,16,38,43,65,70,78\}$ = $R_2$ = $\theta_{81,3,1}(R^{81,7}_1)$ and

\vspace{.1cm}  
$R^{81,7}_3$ = $\{2,3,25,29,52,56,78,79\}$ = $R_3$ = $\theta_{81,3,2}(R^{81,7}_1)$. 

\vspace{.1cm}
$\Rightarrow$ $C_{81}(R_i)$ are isomorphic for $i$ = 1,2,3.

Moreover, 2 is the only positive integer between 1 and 3 and relative prime to 3 and

$2R^{81,7}_1$ = $2 \times \{3,7,20,34,47,61,74,78\}$ = $\{6,14,40,68,13,41,67,75\}$ $\neq$ $R_2, R_3$;

$2R^{81,7}_2$ = $2 \times \{3,11,16,38,43,65,70,78\}$ = $\{6,22,32,76,45,49,59,75\}$ $\neq$ $R_1, R_3$ and 

$2R^{81,7}_3$ = $2 \times \{2,3,25,29,52,56,78,79\}$ = $\{4,6,50,58,23,31,75,77\}$ $\neq$ $R_1, R_2$. 

\vspace{.1cm}
$\Rightarrow$ $C_{81}(R_i)$ are not Type-1 isomorphic, $1 \leq i \leq 3$.

\vspace{.1cm}
$\Rightarrow$ $C_{81}(R_i)$ are Type-2 isomorphic w.r.t. $m$ = 3, $1 \leq i \leq 3$ and

\vspace{.1cm}
$T2_{81,3}(C_{81}(R_1))$ = $\{C_{81}(R_1), C_{81}(R_2), C_{81}(R_3)\}$ = $T2_{81,3}(C_{81}(R_2))$ = $T2_{81,3}(C_{81}(R_3))$. 

Therefore using Theorem \ref{a20}, for $1 \leq i,j \leq 3$, $(T2_{81,3}(C_{81}(R_i)), \circ)$ = $(T2_{81,3}(C_{81}(R_j)), \circ)$ is a subgroup of $(V_{81,3}(C_{81}(R_i)), \circ)$ = $(V_{81,3}(C_{81}(R_j)), \circ)$. 

\begin{figure}[ht]
	\centerline{\includegraphics[width=6.2in]{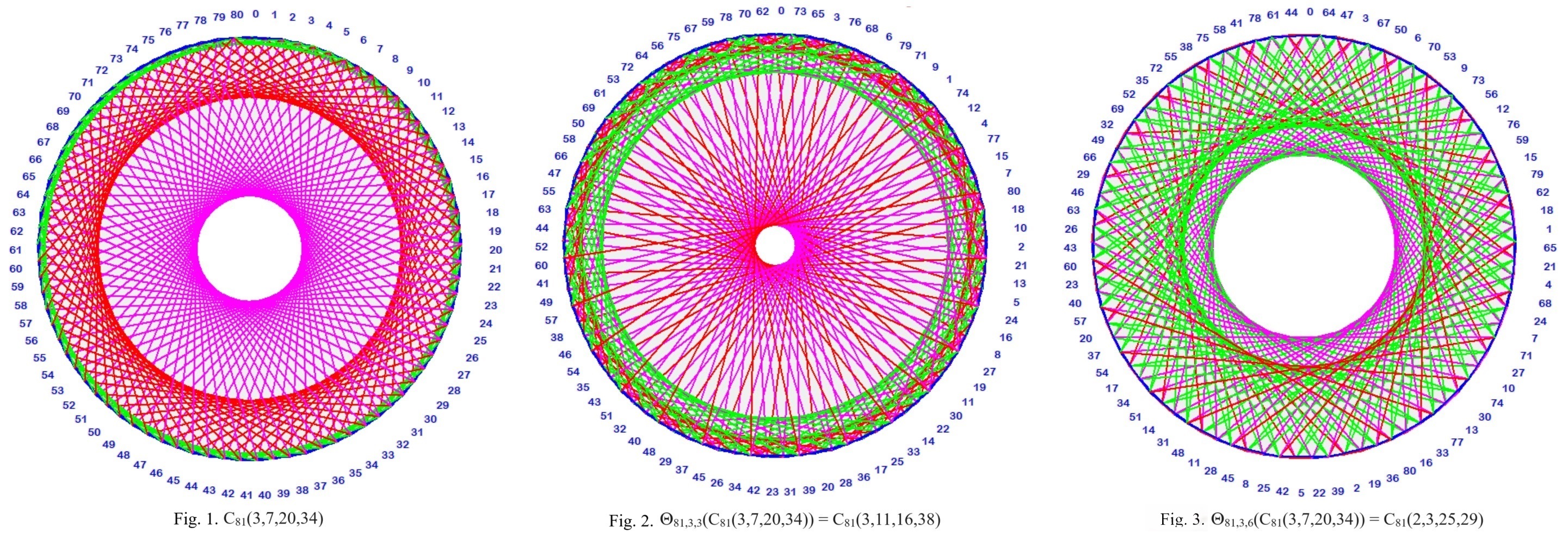}}
\end{figure}
\begin{figure}[ht]
	\centerline{\includegraphics[width=6.2in]{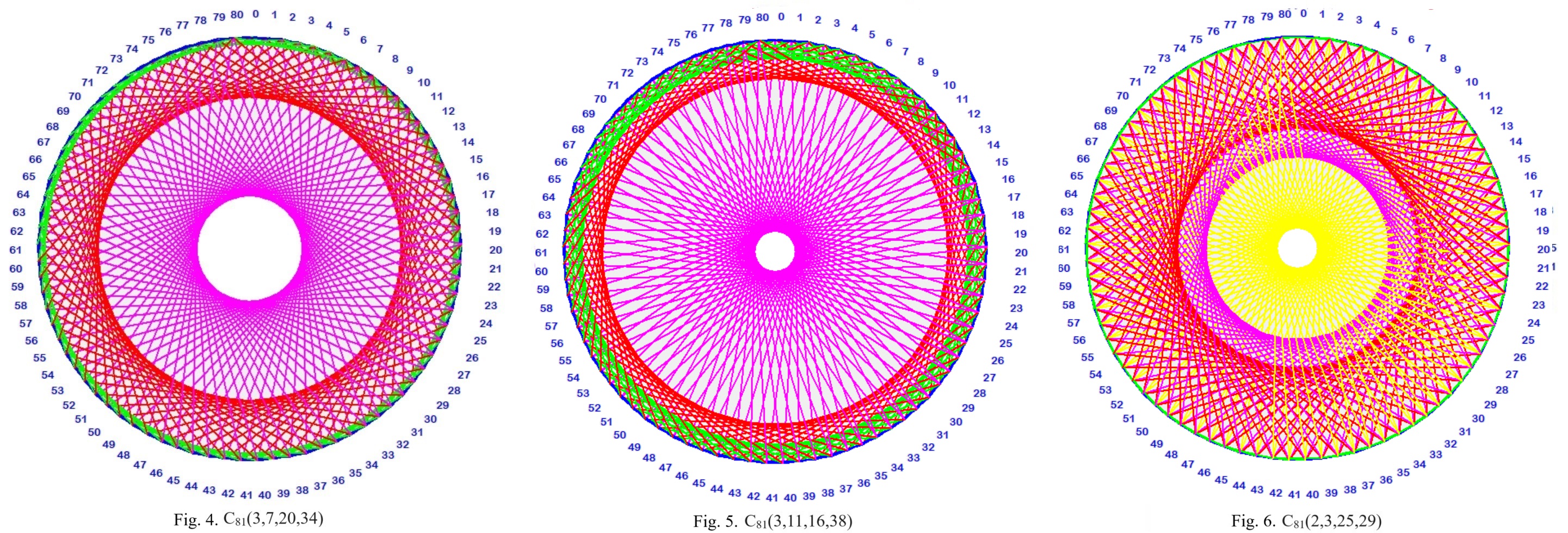}}
\end{figure}
In Figures 1,2,3, isomorphic circulant graphs $\theta_{81,3,0}(C_{81}(R_1))$ = $C_{81}(R_1)$, $\theta_{81,3,3}(C_{81}(R_1))$ = $C_{81}(R_2)$ and $\theta_{81,3,6}(C_{81}(R_1))$ = $C_{81}(R_3)$ of Type-2 w.r.t. $m$ = 3 are given whereas circulant graphs $C_{81}(R_1)$ = $C_{81}(3,7, 20,34)$, $C_{81}(R_2)$ = $C_{81}(3,11, 16,38)$ and $C_{81}(R_3)$ = $C_{81}(2,3, 25,29)$ are given in Figures 4, 5, 6. \hfill $\Box$

\begin{theorem} \label{c4} {\rm Let $p$ be an odd prime number, $1 \leq i \leq p$, $1 \leq x \leq p-1$, $y\in\mathbb{N}_0$, $0 \leq y \leq np - 1$, $1 \leq x+yp \leq np^2-1$, $d^{np^3, x+yp}_i$ = $(i-1)xpn+$ $x+yp$ and $R^{np^3, x+yp}_i$ = $\{p$, $d^{np^3,x+yp}_i$, $np^2-d^{np^3,x+yp}_i$, $np^2+d^{np^3, x+yp}_i$, $2np^2-$ $d^{np^3, x+yp}_i$, $2np^2+d^{np^3, x+yp}_i$, $3np^2-d^{np^3, x+yp}_i$, $3np^2+d^{np^3, x+yp}_i$, . . . , $(p-1)np^2$ - $d^{np^3, x+yp}_i$, $(p-1)np^2+d^{np^3, x+yp}_i$, $np^3-d^{np^3, x+yp}_i$, $np^3-p\}$. Then, for $i$ = 1 to $p$, $T2_{np^3, p}(C_{np^3}(R^{np^3, x+yp}_i))$ = $\{\theta_{np^3,p,jn}(C_{np^3}(R^{np^3,x+yp}_i))$ = $C_{np^3}(R^{np^3,x+yp}_{i+j}) :$ $j$ = $0,1,...,p-1$ and $i+j$ in $C_{np^3}(R^{np^3,x+yp}_{i+j})$ is calculated under addition modulo $p \}$ and $(T2_{np^3, p}(C_{np^3}(R^{np^3, x+yp}_i)), \circ)$ is a Type-2 group of order $p$.}
\end{theorem}
\begin{proof}\quad  Using Theorem \ref{c1}, $\theta_{np^3,p,jn}(C_{np^3}(R^{np^3,x+yp}_i))$ = $C_{np^3}(R^{np^3,x+yp}_{i+j})$ and for $i$ = 1 to $p$, the $p$ circulant graphs $C_{np^3}(R^{np^3,x+yp}_i)$ are isomorphic of Type-2 w.r.t. $p$ when $p$ is an odd prime number, $1 \leq i \leq p$, $1 \leq x \leq p-1$, $y\in\mathbb{N}_0$, $0 \leq y \leq np - 1$, $1 \leq x+yp \leq np^2-1$, $d^{np^3, x+yp}_i$ = $(i-1)xpn+$ $x+yp$, $R^{np^3, x+yp}_i$ = $\{p,$ $d^{np^3,x+yp}_i,$ $np^2-d^{np^3,x+yp}_i,$ $np^2+d^{np^3, x+yp}_i,$ $2np^2-$ $d^{np^3, x+yp}_i,$ $2np^2+d^{np^3, x+yp}_i,$ $3np^2-d^{np^3, x+yp}_i,$ $3np^2+d^{np^3, x+yp}_i,$ ..., $(p-1)np^2$ $-d^{np^3, x+yp}_i,$ $(p-1)np^2+d^{np^3, x+yp}_i,$ $np^3-d^{np^3, x+yp}_i,$ $np^3-p\}$ and $i+j$ in $R^{np^3, x+yp}_{i+j}$ is calculated under addition modulo $p$.
	
This implies, $\theta_{np^3,p,jn}(C_{np^3}(R^{np^3,x+yp}_j))$ = $C_{np^3}(R^{np^3,x+yp}_{i+j})\in T2_{np^3, p}(C_{np^3}(R^{np^3, x+yp}_i))$ for every $i$ and $j$, $1 \leq i,j \leq p$ and $i+j$ in $R^{np^3, x+yp}_{i+j}$ is calculated under addition modulo $p$. Our aim is to prove that $T2_{np^3, p}(C_{np^3}(R^{np^3, x+yp}_i))$ = $\{C_{np^3}(R^{np^3,x+yp}_j): j = 1,2,...,p\}$, $1 \leq i \leq p$. It is enough to prove that $\theta_{np^3,p,t}(C_{np^3}(R^{np^3,x+yp}_j))\notin T2_{np^3, p}(C_{np^3}(R^{np^3, x+yp}_i))$ for $t \neq 0,n,2n,...,(p-1)n$, $0 \leq t \leq np^2-1$, $1 \leq i,j \leq p$. That is to prove that $\theta_{np^3,p,t}(C_{np^3}(R^{np^3,x+yp}_j))\notin T2_{np^3, p}(C_{np^3}(R^{np^3, x+yp}_i))$ for all $t$ = $nqp- s$, $1 \leq q \leq p$, $1 \leq s \leq n-1$ and $1 \leq i,j \leq p$ since $0 \leq t \leq np^2-1$. 

Let $t$ = $nqp - s$, $1 \leq q \leq p$, $1 \leq s \leq n-1$ and $n,s\in\mathbb{N}$. Consider,  
$\theta_{np^3,p,t}(R^{np^3,x+yp}_i)$ = $\theta_{np^3,p,nqp-s}(R^{np^3,x+yp}_i)$ 
 = $\theta_{np^3,p,nqp-s}(\{p, d^{np^3,x+yp}_i, np^2-d^{np^3,x+yp}_i$, $np^2+d^{np^3, x+yp}_i$, 

\hfill $2np^2-d^{np^3, x+yp}_i$, $2np^2+d^{np^3, x+yp}_i$, $3np^2-d^{np^3, x+yp}_i$, $3np^2+d^{np^3, x+yp}_i$, ..., 

\hfill $(p-1)np^2-d^{np^3, x+yp}_i$, $(p-1)np^2+d^{np^3, x+yp}_i$, $np^3-d^{np^3, x+yp}_i$, $np^3-p\})$ 

\hspace{.8cm} = $\{p, np^3-p \}$ $\bigcup$ ($xp(nqp-s) + \{ d^{np^3,x+yp}_i$, $np^2+d^{np^3, x+yp}_i$, $2np^2+d^{np^3, x+yp}_i$, 

\hfill $3np^2+d^{np^3, x+yp}_i$, ..., $(p-1)np^2+d^{np^3, x+yp}_i \}$) 

\hspace{1.2cm} $\bigcup$ ($(p-x)p(nqp-s) + \{ np^2-d^{np^3,x+yp}_i$, $2np^2-d^{np^3, x+yp}_i$, 

\hfill $3np^2-d^{np^3, x+yp}_i$, ..., $(p-1)np^2-d^{np^3, x+yp}_i$, $np^3-d^{np^3, x+yp}_i \}$) 

\hspace{.8cm} = $\{p, np^3-p \}$ $\bigcup$ $\{ xp(nqp-s) +  d^{np^3,x+yp}_i$, $xp(nqp-s) + np^2+d^{np^3, x+yp}_i$,  

\hfill $xp(nqp-s) + 2np^2+d^{np^3, x+yp}_i$, $xp(nqp-s) + 3np^2+d^{np^3, x+yp}_i$, $\dots$, 

\hfill $xp(nqp-s) + (p-1)np^2+d^{np^3, x+yp}_i \}$ 

\hspace{1cm} $\bigcup$ $\{(p-x)p(nqp-s) + np^2-d^{np^3,x+yp}_i$,  $(p-x)p(nqp-s) + 2np^2-d^{np^3, x+yp}_i$, 

\hfill $(p-x)p(nqp-s) + 3np^2-d^{np^3, x+yp}_i$, 
$\dots$, $(p-x)p(nqp-s) + (p-1)np^2-d^{np^3, x+yp}_i$,  

\hfill $(p-x)p(nqp-s) + np^3-d^{np^3, x+yp}_i \}$, $0 \leq x \leq p-1$, $1 \leq i \leq p$, $1 \leq x \leq p-1$, 

\hfill $y\in\mathbb{N}_0$, $0 \leq y \leq np - 1$, $1 \leq x+yp \leq np^2-1$, $d^{np^3, x+yp}_i = (i-1)xpn+$ $x+yp$. 	

\vspace{.2cm}
\noindent
{\it Claim.} For $1 \leq s \leq n-1$, $\theta_{np^3,p,s}(C_{np^3}(R^{np^3,x+yp}_j))\notin T2_{np^3, p}(C_{np^3}(R^{np^3, x+yp}_i))$, $1 \leq i,j \leq p$.

$\theta_{np^3,p,s}(R^{np^3,x+yp}_i)$ = $\theta_{np^3,p,s}(\{p, d^{np^3,x+yp}_i, np^2-d^{np^3,x+yp}_i$, $np^2+d^{np^3, x+yp}_i$, 

\hfill $2np^2-d^{np^3, x+yp}_i$, $2np^2+d^{np^3, x+yp}_i$, $3np^2-d^{np^3, x+yp}_i$, $3np^2+d^{np^3, x+yp}_i$, $\dots$, 

\hfill $(p-1)np^2-d^{np^3, x+yp}_i$, $(p-1)np^2+d^{np^3, x+yp}_i$, $np^3-d^{np^3, x+yp}_i$, $np^3-p\})$ 

\hspace{.8cm} = $\{p, np^3-p \}$ $\bigcup$ ($xps + \{ d^{np^3,x+yp}_i$, $np^2+d^{np^3, x+yp}_i$, $2np^2+d^{np^3, x+yp}_i$, 

\hfill $3np^2+d^{np^3, x+yp}_i$, $\dots$, $(p-1)np^2+d^{np^3, x+yp}_i \}$) 

\hspace{1.2cm} $\bigcup$ ($(p-x)ps + \{ np^2-d^{np^3,x+yp}_i$, $2np^2-d^{np^3, x+yp}_i$, 

\hfill $3np^2-d^{np^3, x+yp}_i$, $\dots$, $(p-1)np^2-d^{np^3, x+yp}_i$, $np^3-d^{np^3, x+yp}_i \}$) 

\hspace{.8cm} = $\{p, np^3-p \}$ $\bigcup$ $\{ xps + d^{np^3,x+yp}_i$, $xps + np^2+d^{np^3, x+yp}_i$, $xps + 2np^2+d^{np^3, x+yp}_i$, 

\hfill $xps + 3np^2+d^{np^3, x+yp}_i$, $\dots$, $xps + (p-1)np^2+d^{np^3, x+yp}_i \}$ 

\hfill $\bigcup$ $\{ (p-x)ps + np^2-d^{np^3,x+yp}_i$, $(p-x)ps + 2np^2-d^{np^3, x+yp}_i$, $(p-x)ps + 3np^2-d^{np^3, x+yp}_i$,

\hfill  $\dots$, $(p-x)ps + (p-1)np^2-d^{np^3, x+yp}_i$, $(p-x)ps + np^3-d^{np^3, x+yp}_i \}$, $1 \leq i \leq p$,

\hfill  $1 \leq s,x \leq p-1$,  $y\in\mathbb{N}_0$, $0 \leq y \leq np - 1$, $1 \leq x+yp \leq np^2-1$, $d^{np^3, x+yp}_i = (i-1)xpn+$ $x+yp$. 	

\vspace{.2cm}
\noindent
{\it Sub-claim.} For $s$ = 1, $\theta_{np^3,p,s}(C_{np^3}(R^{np^3,x+yp}_j))\notin T2_{np^3, p}(C_{np^3}(R^{np^3, x+yp}_i))$, $1 \leq i,j \leq p$.

$\theta_{np^3,p,1}(R^{np^3,x+yp}_i)$ = $\{p, np^3-p \}$ $\bigcup$ $\{ xp + d^{np^3,x+yp}_i$, $xp + np^2+d^{np^3, x+yp}_i$, $xp + 2np^2+d^{np^3, x+yp}_i$, 

\hfill $xp + 3np^2+d^{np^3, x+yp}_i$, $\dots$, $xp + (p-1)np^2+d^{np^3, x+yp}_i \}$ 

\hfill $\bigcup$ $\{ (p-x)p + np^2-d^{np^3,x+yp}_i$, $(p-x)p + 2np^2-d^{np^3, x+yp}_i$, $(p-x)p + 3np^2-d^{np^3, x+yp}_i$,

\hfill  $\dots$, $(p-x)p + (p-1)np^2-d^{np^3, x+yp}_i$, $(p-x)p + np^3-d^{np^3, x+yp}_i \}$, $1 \leq i \leq p$,

  $1 \leq x \leq p-1$,  $y\in\mathbb{N}_0$, $0 \leq y \leq np - 1$, $1 \leq x+yp \leq np^2-1$, $d^{np^3, x+yp}_i $= $(i-1)xpn+$ $x+yp$. 	

  Let $z\in\mathbb{N}$ $\ni$ $1 \leq z \leq p$. Among the elements of the above set, 
 \\
  $(xp +  d^{np^3,x+yp}_i)$ + $((p-x)p + (p-z)np^2-d^{np^3, x+yp}_i)$ = $p^2 + (p-z)np^2$ = $np^3 + (1-nz)p^2$ $\neq$ $np^3$ for any $z\in\mathbb{N}$ $\ni$ $1 \leq z \leq p$. This implies, the elements of the set $\theta_{np^3,p,1}(R^{np^3,x+yp}_i)$ are not satisfying the symmetric equidistance condition. Therefore, using Theorem \ref{ab14}, $\theta_{np^3,p,1}(R^{np^3,x+yp}_i) \neq C_{np^3}(R)$ for any $R \subseteq [1, np^3/2]$ and thereby $\theta_{np^3,p,1}(R^{np^3,x+yp}_i) \notin T2_{np^3,p}(R^{np^3,x+yp}_i)$, $1 \leq i \leq p$. Hence the sub-claim is true.
  
  In a similar way, when  $1 \leq s \leq n-1$ and $n,s\in\mathbb{N}$, we have 
  
  $\theta_{np^3,p,s}(R^{np^3,x+yp}_i)$ = $\{p, np^3-p \}$ $\bigcup$ $\{ xps + d^{np^3,x+yp}_i$, $xps + np^2+d^{np^3, x+yp}_i$,  
    
  \hfill $xps + 2np^2+d^{np^3, x+yp}_i$, $xps + 3np^2+d^{np^3, x+yp}_i$, $\dots$, $xps + (p-1)np^2+d^{np^3, x+yp}_i \}$ 
  
  \hfill $\bigcup$ $\{ (p-x)ps + np^2-d^{np^3,x+yp}_i$, $(p-x)ps + 2np^2-d^{np^3, x+yp}_i$, $(p-x)ps + 3np^2-d^{np^3, x+yp}_i$,
  
  \hfill  $\dots$, $(p-x)ps + (p-1)np^2-d^{np^3, x+yp}_i$, $(p-x)ps + np^3-d^{np^3, x+yp}_i \}$, $1 \leq i \leq p$,
  
 $1 \leq s,x \leq p-1$,  $y\in\mathbb{N}_0$, $0 \leq y \leq np - 1$, $1 \leq x+yp \leq np^2-1$, $d^{np^3, x+yp}_i = (i-1)xpn+$ $x+yp$. 	
    
 For any $z\in\mathbb{N}$ $\ni$ $1 \leq z \leq p$, among the elements of the above set, 
\\
$(xps +  d^{np^3,x+yp}_i)$ + $((p-x)ps + (p-z)np^2-d^{np^3, x+yp}_i)$ = $sp^2 + (p-z)np^2$ $\neq$ $np^3$, $1 \leq s \leq p-1$. This implies, the elements of the set $\theta_{np^3,p,s}(R^{np^3,x+yp}_i)$ are not satisfying the symmetric equidistance condition, $1 \leq s \leq n-1$. Therefore, using Theorem \ref{ab14}, $\theta_{np^3,p,s}(R^{np^3,x+yp}_i) \notin T2_{np^3,p}(R^{np^3,x+yp}_i)$ when $1 \leq s \leq n-1$ and $1 \leq i \leq p$. Hence the claim is true in this case.

Similarly, we can show that $\theta_{np^3,p,nqp-s}(R^{np^3,x+yp}_i) \notin T2_{np^3,p}(R^{np^3,x+yp}_i)$ when $1 \leq s \leq n-1$, $1 \leq q \leq p$ and $n,p,q,s\in\mathbb{N}$.

This implies that for $i$ = 1, 2, . . . , $p$, $T2_{np^3, p}(C_{np^3}(R^{np^3, x+yp}_i))$ = $\{\theta_{np^3,p,jn}(C_{np^3}(R^{np^3,x+yp}_i))$ = $C_{np^3}(R^{np^3,x+yp}_{i+j}) : j = 0,1,...,p-1$ and $i+j$ in $C_{np^3}(R^{np^3,x+yp}_{i+j})$ is calculated under addition modulo $p \}$ and $(T2_{np^3, p}(C_{np^3}(R^{np^3, x+yp}_i)), \circ)$ is a Type-2 group of order $p$. 
\end{proof}

Let us see the importance of the above theorems by the following problem. 

\begin{prm} \label{c5} {\rm Check whether the circulant graphs $C_{1715}(G_i)$ for $i$ = 1 to 7 are isomorphic or not? And classify, if they are isomorphic where 

$G_1$ = $\{7,17,228,262,473,507,718,752\}$,

$G_2$ = $\{7,122,123,367,368,612,613,857\}$,  

$G_3$ = $\{7,18,227,263,472,508,717,753\}$,

$G_4$ = $\{7,87,158,332,403,577,648,822\}$,

$G_5$ = $\{7,53,192,298,437,543,682,788\}$,

$G_6$ = $\{7,52,193,297,438,542,683,787\}$ and

$G_7$ = $\{7,88,157,333,402,578,647,823\}$.}  
\end{prm}
\noindent
{\bf Solution} All the 7 sets $G_1$, $G_2$, . . . , $G_7$ are subsets of $\mathbb{Z}_{1715/2}$ and have 7 as a common element and thereby there is a possibility that $C_{1715}(G_i)$ may be of Type-2 isomorphic w.r.t.  $m$ = 7.  Using Theorem \ref{c1}, we show that circulant graphs $C_{1715}(G_i)$ are Type-2 isomorphic w.r.t. $m$ = $p$ = 7,  $1 \leq i \leq p$ = 7.  

Here, 1715 = $5\times {7}^3$. Let $n = 5$ and $p = 7$. This implies, $np$ = 35, $np^2 = 245$ and $np^3$ = 1715. 

The minimum jump size, other than 7, in the seven $G_i$ sets is 17 which implies, as in Theorem \ref{c1}, 17 = $x+yp$ = $3 + 2\times 7$. This implies, $x = 3$, $y = 2$ and 

$d^{np^3,x+yp}_i$ = $d^{5\times {7}^3,3+2\times 7}_i$ = $(i-1)xpn+yp+x$ = $105(i-1)+17$, $i = $ 1 to 7. This implies, 

$d^{5\times {7}^3,3+2\times 7}_1$ = $17\in G_1$; $d^{5\times {7}^3,3+2\times 7}_2$ = $122\in G_2$; $d^{5\times {7}^3,3+2\times 7}_3$ = $227\in G_3$; 

$d^{5\times {7}^3,3+2\times 7}_4$ = $332\in G_4$; $d^{5\times {7}^3,3+2\times 7}_5$ = $437\in G_5$; $d^{5\times {7}^3,3+2\times 7}_6$ = $542\in G_6$ and

$d^{5\times {7}^3,3+2\times 7}_7$ = $647\in G_7$.
\\
For $i =$ 1 to $p$, we have 

$R^{np^3,yp+x}_i$ = $\{p,$ $d^{np^3,yp+x}_i$, $np^2-d^{np^3,yp+x}_i$, $np^2+d^{np^3,yp+x}_i$, $2np^2-d^{np^3,yp+x}_i$, $2np^2+d^{np^3,yp+x}_i$,

\hfill . . . , $(p-1)np^2-d^{np^3,yp+x}_i$, $(p-1)np^2+d^{np^3,yp+x}_i$, $np^3-d^{np^3,yp+x}_i$, $np^3-p\}$. 
\\
This implies, under modulo 1715,

$R^{5\times 7^3,2\times 7+3}_1$ = $\{7, 17, 228, 262, 473, 507, 718, 752, $ 

\hfill $963, 997, 1208, 1242, 1453, 1487, 1698, 1708\}$ = $G_1 \cup (1715 - G_1)$; 

$R^{5\times 7^3,2\times 7+3}_2$ = $\{7, 122, 123, 367, 368, 612, 613, 857, $ 

\hfill $858, 1102, 1103, 1347, 1348, 1592, 1593, 1708\}$ = $G_2 \cup (1715 - G_2)$; 

$R^{5\times 7^3,2\times 7+3}_3$ = $\{7, 227, 18, 472, 263, 717, 508, 962, 753, 1207, 998,$

\hfill $1452, 1243, 1697, 1488, 1708\}$ = $G_3 \cup (1715 - G_3)$; 

$R^{5\times 7^3,2\times 7+3}_4$ = $\{7, 332, (245-332)+1715 = 1628, 577, 158, 822, 403, 1067,$ 

\hfill $648, 1312, 893, 1557, 1138,  1802 = 87, 1383, 1708\}$ = $G_4 \cup (1715 - G_4)$; 

$R^{5\times 7^3,2\times 7+3}_5$ = $\{7, 437, (245-437)+1715 = 1523, 682, 53, 927,$ 

\hfill  $298, 1172, 543, 1417, 788, 1662, 1033, 192, 1278, 1708\}$ = $G_5 \cup (1715 - G_5)$; 

$R^{5\times 7^3,2\times 7+3}_6$ = $\{7, 542, (245-542)+1715 = 1418, 787, 1663, 1032,$ 

\hfill $ 193, 1277, 438, 1522, 683, 52, 928, 297, 1173, 1708\}$ = $G_6 \cup (1715 - G_6)$;

$R^{5\times 7^3,2\times 7+3}_7$ = $\{7, 647, (245-647)+1715 = 1313, 892, 1558, 1137,$ 

\hfill $88, 1382, 333, 1627, 578, 157, 823, 402, 1068, 647, 1708\}$ = $G_7 \cup (1715 - G_7)$. 

\vspace{.2cm}
\noindent
{\bf{\it Claim 1.}} $\theta_{1715,7,n}(R^{1715,17}_1)$ = $R^{1715,17}_2$, $n$ = 5.   

\vspace{.2cm}
\noindent
$\theta_{1715,7,n}(R^{1715,17}_1)$ = $\theta_{1715,7,5}(R^{1715,17}_1)$ = $\theta_{1715,7,5}(\{7, 17, 228, 262, 473,$ 

~\hfill $507, 718, 752, 963, 997, 1208, 1242, 1453, 1487, 1698, 1708 \})$ 

 = $\theta_{1715,7,5}(\{7, 1708\})$ $\bigcup$ $\theta_{1715,7,5}(\{17, 262, 507, 752, 997, 1242, 1487\})$ 

~\hfill $\bigcup$ $\theta_{1715,7,5}(\{228, 473, 718, 963, 1208, 1453, 1698\})$

= $\{7, 1708\} \bigcup (3\times 7 \times 5+\{17, 262, 507, 752, 997, 1242, 1487\})$

~\hfill $\bigcup$ $(4\times 7 \times 5+\{228, 473, 718, 963, 1208, 1453, 1698\})$ 

 = $\{7, 1708\}$ $\bigcup$ $\{122, 367, 612, 857, 1102, 1347, 1592\} 
\bigcup \{368, 613, 858, 1103, 1348, 1593, 123\}$ = $R^{1715,17}_2$.

Hence Claim 1 is true.

\vspace{.2cm}
\noindent
{\it Claim 2.} $\theta_{1715,7,4n}(R^{1715,17}_1)$ = $R^{1715,17}_5$, $n = 5$.   

\vspace{.2cm}
\noindent
 $\theta_{1715,7,4n}(R^{1715,17}_1)$ = $\theta_{1715,7,20}(\{7, 1708\}) \bigcup \theta_{1715,7,20}(\{17, 262, 507, 752, 997, 1242, 1487\})$ 

~\hfill $\bigcup$ $\theta_{1715,7,20}(\{228, 473, 718, 963, 1208, 1453, 1698\})$

 = $\{7, 1708\}$ $\bigcup$ $(3\times 20 \times 7+\{17, 262, 507, 752, 997, 1242, 1487\})$

~\hfill $\bigcup$ $(4\times 20 \times 7+\{228, 473, 718, 963, 1208, 1453, 1698\})$ 

 = $\{7, 1708\}$ $\bigcup$ $\{437, 682, 927, 1172, 1417, 1662, 1907 = 192\}$ 

~\hfill $\bigcup$ $\{788, 1033, 1278, 1523, 1768 = 53, 2013 = 298, 543\}$ = $R^{1715,17}_5$.

Hence Claim 2 is true.

$\Rightarrow$ $\theta_{1715,7,n}(C_{1715}(G_1))$ = $\theta_{1715,7,n}(C_{1715}(R^{1715,17}_1))$ = $C_{1715}(R^{1715,17}_2)$ = $C_{1715}(G_2)$ and 

\hspace{.25cm} $\theta_{1715,7,4n}(C_{1715}(G_1))$ = $\theta_{1715,7,4n}(C_{1715}(R^{1715,17}_1))$ = $C_{1715}(R^{1715,17}_5)$ = $C_{1715}(G_5)$.

$\Rightarrow$ $C_{1715}(G_1)$, $C_{1715}(G_2)$ and $C_{1715}(G_5)$ are isomorphic graphs.

Similarly, we can show that $\theta_{1715,7,jn}(C_{1715}(G_1))$ = $C_{1715}(G_{j+1})$ as well as $\theta_{1715,7,jn}(C_{1715}(G_i))$ = $C_{1715}(G_{i+j})$ where $i+j$ in $G_{i+j}$ is calculated under addition modulo 7, $0 \leq i,j \leq 6$ and $G_0$ = $G_7$. This implies, $C_{1715}(G_i)$ are isomorphic circulant graphs for $i$ = 1 to 7.

Also, for $1 \leq i,j \leq 7$, $i \neq j$, $1 < s < 7$ and $s\in\varphi_{7}$, $sR^{1715,17}_i$ $\neq$ $R^{1715,17}_i, R^{1715,17}_j$ since $7\in R^{1715,17}_i$ but $7\notin sR^{1715,17}_i$ whereas $7s\in sR^{1715,17}_i$ but $7s \notin R^{1715,17}_i$. This implies, $C_{1715}(R^{1715,17}_i)$ and $C_{1715}(R^{1715,17}_j)$ are not of Type-1 isomorphic when $i \neq j$ and $1 \leq i,j \leq 7$. This implies, $C_{1715}(G_i)$ and $C_{1715}(G_{j})$ are isomorphic but not of Type-1 isomorphic when $i \neq j$ and $1 \leq i,j \leq 7$. 

This implies, $\theta_{1715,7,jn}(C_{1715}(G_1))$ = $\theta_{1715,7,jn}(C_{1715}(R^{1715,17}_1))$ = $C_{1715}(R^{1715,17}_{j+1})$ = $C_{1715}(G_{j+1})$ are isomorphic circulant graphs of Type-2 w.r.t.  $m$ = 7, $0 \leq j \leq p-1$ = 6.  \hfill $\Box$

\begin{lemma} \label{c6} {\rm Let $p$ be an odd prime, $n\in\mathbb{N}$, $k \geq 3$, $1 \leq x \leq p - 1$, $1 \leq x+yp \leq np^2 - 1$, $0 \leq y \leq np-1$, $1 \leq i \leq p$, $d^{np^3,x+yp}_i$ = $(i-1)xpn+x+yp$ and $R^{np^3,x+yp}_i$ = $\{d^{np^3,x+yp}_i,$ $np^2-d^{np^3,x+yp}_i$, $np^2+$ $d^{np^3,x+yp}_i$, $2np^2-d^{np^3,x+yp}_i,$ $2np^2+$ $d^{np^3,x+yp}_i,$ $3np^2-d^{np^3,x+yp}_i,$ $3np^2+$ $d^{np^3,x+yp}_i,$ $\dots$, $(p-1)np^2-d^{np^3,x+yp}_i,$ $(p-1)np^2+d^{np^3,x+yp}_i,$ $np^3-d^{np^3,x+yp}_i,$ $pp_1,pp_2,...,pp_{k-2},$ $p(np^3-p_{k-2}),$ $p(np^3-p_{k-3}),$ . . . , $p(np^3-p_1)\}$ where $\gcd(p_1,p_2,...,p_{k-2}) = 1$, $y\in\mathbb{N}_0$ and $i,k,x,p_1,p_2,...,p_{k-2} \in \mathbb{N}.$ Then, for $1 \leq i \leq p$ and a given set of values of $k$, $n$, $p$, $x$, $y$, $p_1,$ $p_2,$ . . . , $p_{k-3}$ and $p_{k-2}$, $R^{np^3,np^2-x-yp}_i$ = $R^{np^3,x+yp}_i$.} 
\end{lemma} 

\begin{proof} For $1 \leq x \leq p - 1$, $1 \leq x+yp \leq np^2 - 1$, $0 \leq y \leq np - 1$, $1 \leq i \leq p$ and $y\in\mathbb{N}_0$, we have $d^{np^3,x+yp}_i = (i-1)xpn+x+yp$ and $d^{np^3,np^2-x-yp}_i = (i-1)(p-x)pn+np^2-x-yp$ = $inp^2-d^{np^3,x+yp}_i$. 
	
Therefore $R^{np^3,np^2-x-yp}_i$ = $\{inp^2-d^{np^3,x+yp}_i,$ $(p-(i-1))np^2+d^{np^3,x+yp}_i$, $(i+1)np^2-d^{np^3,x+yp}_i$, $(p-(i-2))np^2+d^{np^3,x+yp}_i,$ $(i+2)np^2-d^{np^3,x+yp}_i$, $(p-(i-3))np^2+d^{np^3,x+yp}_i,$ $(i+3)np^2-d^{np^3,x+yp}_i,$ ..., $(p-(i+1))np^2+d^{np^3,x+yp}_i,$ $(i+p-1)np^2-d^{np^3,x+yp}_i,$ $(p-i)np^2+d^{np^3,x+yp}_i,$ $pp_1,$ $pp_2,$ . . . , $pp_{k-2},$ $p(np^3-p_{k-2}),$ $p(np^3-p_{k-3}),$ . . . , $p(np^3-p_1)\}$ = $R^{np^3,x+yp}_i$ since $\{p-(i-1)$, $p-(i-2)$, $p-(i-3)$, . . . , $p-(i+1) = p-i+p-1$, $p-i = p-i+p\}$ = $\{i$, $i+1$, $i+2$, . . . , $i+p-2$, $i+p-1\}$ = $\{1, 2, 3, . . . , p-1, p = 0\}$ under addition modulo $p$, $1 \leq i \leq p$. Hence the lemma is proved.
\end{proof}

\par
We get the following as a particular case of the above lemma when $n$ = 1.

\begin{lemma} \label{c7} {\rm Let $p$ be an odd prime, $1 \leq x \leq p - 1$, $0 \leq y \leq p - 1$, $1 \leq x+yp \leq p^2 - 1$, $k \geq 3$, $1 \leq i \leq p$, $d^{p^3,x+yp}_i = (i-1)xp+x+yp$ and $R^{p^3,x+yp}_i$ = $\{d^{p^3,x+yp}_i,$ $p^2-d^{p^3,x+yp}_i$, $p^2+d^{p^3,x+yp}_i$, $2p^2-d^{p^3,x+yp}_i,$ $2p^2+d^{p^3,x+yp}_i,$ $3p^2-d^{p^3,x+yp}_i,$ $3p^2+d^{p^3,x+yp}_i,$ . . . , $(p-1)p^2-d^{p^3,x+yp}_i,$ $(p-1)p^2+d^{p^3,x+yp}_i,$ $p^3-d^{p^3,x+yp}_i,$ $pp_1,$ $pp_2,$ . . . , $pp_{k-2},$ $p(p^3-p_{k-2}),$ $p(p^3-p_{k-3}),$ . . . , $p(p^3-p_1)\}$ where $\gcd(p_1,p_2,...,p_{k-2}) = 1$, $y\in\mathbb{N}_0$ and $i,k,x,p_1,p_2,...,p_{k-2} \in \mathbb{N}.$ Then, for a given set of values of $k, p, x, y, p_1, p_2,  . . . , p_{k-3}$ and $p_{k-2}$, $R^{p^3,p^2-x-yp}_i$ = $R^{p^3,x+yp}_i$, $1 \leq i \leq p$. \hfill $\Box$}
\end{lemma} 

\par
In the above lemma when $y$ = 0, we get the following.

\begin{lemma} \label{c8} {\rm Let $p$ be an odd prime, $1 \leq x \leq p - 1$, $k \geq 3$, $1 \leq i \leq p$, $d^{p^3, x}_i = (i-1)xp+x$ and $R^{p^3,x}_i$ = $\{d^{p^3, x}_i,$ $p^2-d^{p^3, x}_i$, $p^2+d^{p^3, x}_i$, $2p^2-d^{p^3, x}_i,$ $2p^2+d^{p^3, x}_i,$ $3p^2-d^{p^3, x}_i,$ $3p^2+d^{p^3, x}_i,$ . . . , $(p-1)p^2-d^{p^3, x}_i,$ $(p-1)p^2+d^{p^3, x}_i,$ $p^3-d^{p^3, x}_i,$ $pp_1,$ $pp_2,$ . . . , $pp_{k-2},$ $p(p^3-p_{k-2}),$ $p(p^3-p_{k-3}),$ . . . , $p(p^3-p_1)\}$ where $\gcd(p_1,p_2,...,p_{k-2}) = 1$ and $i,k,x,p_1,p_2,...,p_{k-2} \in \mathbb{N}.$ Then, for a given set of values of $k, p, x, p_1, p_2,  . . . , p_{k-3}$ and $p_{k-2}$, $R^{p^3,p-x}_i$ = $R^{p^3,x}_i$, $1 \leq x \leq p-1$ and $1 \leq i \leq p$. \hfill $\Box$}
\end{lemma} 

\begin{note}\label{c9} \quad {\rm If we put $p$ = 2 and consider $n \geq 2$ in Theorems \ref{c1} and \ref{c4}, then $x$ = 1, $0 \leq y \leq 2n-1$, $n \geq 2$, $y\in\mathbb{N}_0$ and the above results become Theorems \ref{a17c} and \ref{a20}, respectively and thereby we get the corresponding more general result as follows. Also, using remark \ref{r2.6} in Theorems \ref{c10} and \ref{c12}, we consider $pq$ in the place of jump size $p$ in the circulant graphs where $q\in\mathbb{N}$ and $\gcd(np^3, pq)$ = $p$ = $m$. }
\end{note}

\begin{theorem} \label{c10} {\rm Let $p$ be a prime number, $1 \leq x \leq p-1$, $1 \leq i \leq p$, $y\in\mathbb{N}_0$, $0 \leq y \leq np-1$, $1 \leq x+yp \leq np^2-1$, $d^{np^3, x+yp}_i$ = $(i-1)xpn+x+yp$, $R^{np^3, x+yp}_i$ = $\{pq$, $d^{np^3, x+yp}_i$, $np^2-d^{np^3, x+yp}_i$, $np^2+d^{np^3, x+yp}_i$, $2np^2-d^{np^3, x+yp}_i$, $2np^2+$ $d^{np^3, x+yp}_i,$ $3np^2-d^{np^3, x+yp}_i$, $3np^2+d^{np^3, x+yp}_i$, . . . , $(p-1)np^2$ - $d^{np^3, x+yp}_i$, $(p-1)np^2+d^{np^3, x+yp}_i$, $np^3-d^{np^3, x+yp}_i$, $np^3-pq\}$, $m$ = $\gcd(np^3, pq)$ = $p$, $i,n,q,x\in\mathbb{N}$ and $n \geq 2$ when $p$ = 2. Then, for a given set of values of $n, p, q, x$ and $y$, $(i)$ $\theta_{np^3,p,jn} (C_{np^3}(R^{np^3, x+yp}_i))$ = $C_{np^3}(R^{np^3, x+yp}_{i+j})$ and $(ii)$ the $p$ circulant graphs $C_{np^3}(R^{np^3, x+yp}_i)$ are Type-2  isomorphic w.r.t.  $m$ = $p$ (and without CI-property) where $i+j$ in $R^{np^3, x+yp}_{i+j}$ is calculated under addition modulo $p$, $1 \leq i,j \leq p$. 

\hfill $\Box$}
\end{theorem}

\begin{theorem} \label{c11} {\rm Let $p$ be a prime number, $k \geq 3$, $1 \leq i,j \leq p$, $1 \leq x \leq p-1$, $1 \leq x+yp \leq np^2-1$, $y\in\mathbb{N}_0$, $0 \leq y \leq np - 1$,  $d^{np^3, x+yp}_i$ = $(i-1)xpn+$ $x+yp$, $R^{np^3, x+yp}_i$ = $\{d^{np^3, x+yp}_i,$ $np^2-d^{np^3, x+yp}_i$, $np^2+d^{np^3, x+yp}_i$, $2np^2-$ $d^{np^3, x+yp}_i$, $2np^2+d^{np^3, x+yp}_i$, $3np^2-d^{np^3, x+yp}_i$, $3np^2+d^{np^3, x+yp}_i$, ..., $(p-1)np^2-d^{np^3, x+yp}_i$, $(p-1)np^2+d^{np^3, x+yp}_i$, $np^3-d^{np^3, x+yp}_i$, $pp_1$, $pp_2$, $\dots$, $pp_{k-2}$, $p(np^3-p_{k-2}), p(np^3-p_{k-3}), . . . ,$ $p(np^3-p_1)\}$, $\gcd(p_1,p_2,...,p_{k-2})$ = 1, $i,j,k,n,x,p_1,p_2,...,p_{k-2} \in \mathbb{N}$ and $n \geq 2$ when $p$ = 2. Then, corresponding to each value of $k \geq 3$ and for a given set of values of $p, x, y, p_1, p_2, \dots, p_{k-3}, p_{k-2}, n$
and for $i$ = $1, 2, \dots, p$, $(i)$ $\theta_{np^3,p,jn}(C_{np^3}(R^{np^3, x+yp}_i))$ = $C_{np^3}(R^{np^3, x+yp}_{i+j})$ and $(ii)$ the $p$ circulant graphs $C_{np^3}(R^{np^3, x+yp}_i)$ are either all Type-1 isomorphic or all Type-2 isomorphic w.r.t. $p$ (and without CI-property) where $i+j$ in $R^{np^3, x+yp}_{i+j}$ is calculated under addition modulo $p$, $1 \leq i,j \leq p$.  \hfill $\Box$} 
\end{theorem}

\begin{theorem} \label{c12} {\rm Let $p$ be a prime number, $1 \leq x \leq p-1$, $1 \leq i \leq p$, $y\in\mathbb{N}_0$, $0 \leq y \leq np-1$, $1 \leq x+yp \leq np^2-1$, $d^{np^3, x+yp}_i = (i-1)xpn+$ $x+yp$, $R^{np^3, x+yp}_i$ = $\{pq$, $d^{np^3,x+yp}_i,$ $np^2 - d^{np^3,x+yp}_i,$ $np^2+d^{np^3, x+yp}_i,$ $2np^2-$ $d^{np^3, x+yp}_i,$ $2np^2+d^{np^3, x+yp}_i,$ $3np^2-d^{np^3, x+yp}_i,$ $3np^2+d^{np^3, x+yp}_i,$ ..., $(p-1)np^2$ $-d^{np^3, x+yp}_i,$ $(p-1)np^2+d^{np^3, x+yp}_i,$ $np^3-d^{np^3, x+yp}_i,$ $np^3-pq\}$, $m$ = $\gcd(np^3, pq)$ = $p$, $i,n,q,x\in \mathbb{N}$ and $n \geq 2$ when $p$ = 2. Then, for $i$ = 1 to $p$, $T2_{np^3, m}(C_{np^3}(R^{np^3, x+yp}_i))$ = $\{\theta_{np^3,p,jn}(C_{np^3}(R^{np^3,x+yp}_i))$ = $C_{np^3}(R^{np^3,x+yp}_{i+j}):$ $j$ = $0,1,\dots,p-1$ and $i+j$ in $C_{np^3}(R^{np^3,x+yp}_{i+j})$ is calculated under addition modulo $p\}$ and $(T2_{np^3, m}(C_{np^3}(R^{np^3, x+yp}_i)), \circ)$ is a Type-2 group of order $p$. \hfill $\Box$}
\end{theorem}

\section {Conclusion}

In this paper, we proved that for $i = 1,2,...,p$ and for a given set of values of $n, p, q, x$ and $y$, the $p$ circulant graphs $C_{np^3}(R^{np^3,x+yp}_i)$ are isomorphic of Type-2 w.r.t. $m$ = $p$ and they form the Abelian subgroup $(T2_{np^3,m}(C_{np^3}(R^{np^3,x+yp}_j)),~ \circ)$ of $(V_{np^3,m}(C_{np^3}(R^{np^3,x+yp}_j)),~ \circ)$ where $p$ is a prime number, $d^{np^3, x+yp}_i$ = $(i-1)xpn+x+yp$,  $R^{np^3, x+yp}_i$ = $\{pq$, $d^{np^3, x+yp}_i$, $np^2-d^{np^3, x+yp}_i$, $np^2+d^{np^3, x+yp}_i$, $2np^2-d^{np^3, x+yp}_i$, $2np^2+$ $d^{np^3, x+yp}_i,$ $3np^2-d^{np^3, x+yp}_i$, $3np^2+d^{np^3, x+yp}_i,$ . . . , $(p-1)np^2$ - $d^{np^3, x+yp}_i$, $(p-1)np^2+d^{np^3, x+yp}_i,$ $np^3-d^{np^3, x+yp}_i,$ $np^3-pq\}$, $1 \leq x \leq p-1$, $1 \leq i \leq p$, $y\in\mathbb{N}_0$, $0 \leq y \leq np-1$, $1 \leq x+yp \leq np^2-1$, $i,j,n,q,x\in\mathbb{N}$, $m$ = $\gcd(np^3, pq)$ = $p$ and $n \geq 2$ when $p$ = 2. Type-2 isomorphic circulant graphs are without CI-property. Using Theorem \ref{c4}, a list of Type-2 sets  $T2_{np^3,p}(C_{np^3}(R^{np^3,x+yp}_i))$, each set of order $p$, of isomorphic circulant graphs is given in the Annexure for $p$ = 3,5,7, $n$ = 1,2 and $y$ = 0 whereas in \cite{v24}, more such Type-2 sets are given for $p$ = 3,5,7,11 and $n$ = 1 to 5 and also for $p$ = 13 and $n$ = 1 to 3 where $1 \leq i \leq p$, $1 \leq x \leq p-1$, $1 \leq x+yp \leq np^2 - 1$, $0 \leq y \leq np-1$, $y\in\mathbb{N}_0$ and $p,np^3-p\in R^{np^3,x+yp}_i$. The importance of this theory is that without this theory, it is not easy to establish isomorphism among these circulant graphs when their orders are large. The authors feel that a lot of scope is there to extend many results of this paper. 

\vspace{.1cm}
\noindent
\textbf{Declaration of competing interest}\quad 
The authors declare that they have no conflict of interest.

\vspace{.5cm}
\noindent
\textbf{Acknowledgements}\quad We express our sincere thanks to Prof. S. Krishnan (late), Prof. V. Mohan and Prof. R. Aravamuthan (late), Thiyagarayar College of Engineering, Madurai, Tamil Nadu, India; Prof. M. I. Jinnah (late), Mr. Albert Joseph and Prof. L. John, University of Kerala, Trivandrum, Kerala, India; Prof. K. Vareethiah Konstantine and Prof. S. Amirthaiyan, St. Jude's College, Thoothoor, K. K. District, Tamil Nadu; Mr. R. Benadict, Headmaster (Rtd), Pius XI Higher Secondary School, Thoothoor; Dr. Oscar Fredy, Royal Liverpool University Hospital, Liverpool, U.K. and Prof. Lowell W Beineke, Purdue University, USA for their help and guidance and our sincere thanks to Dr. A. Christopher and Mr. R. Satheesh of S.T. Hindu College, Nagercoil, India for their assistance to develop the VB programs. We also express our gratitude to the Central University of Kerala, Kasaragod, Kerala; St. Jude's College, Thoothoor and S. T. Hindu College, Nagercoil; and Lerroy Wilson Foundation, India (www.WillFoundation.co.in) for providing facilities to do this research work.


\newpage
\begin{center}
	\textbf{ANNEXURE}
\end{center}

\textbf{List of families of Type-2 isomorphic circulant graphs} 

\vspace{.2cm}
$T2_{np^3,p}(C_{np^3}(R^{np^3,x+yp}_i)) = \{C_{np^3}(R^{np^3,x+yp}_i):~i = 1,2,...,p\}$, 

~\hfill $1 \leq x \leq p-1$, $0 \leq y \leq np-1$, $n,x\in\mathbb{N}$, $y\in\mathbb{N}_0$ and $p,np^3\in R^{np^3,x+yp}_i$.

-----------------------------------------------------------------------------

\vspace{.1cm}
$T2_{27,3}(C_{27}(R^{27,1}_i))$, $p$ = 3, $x$ = 1, $y$ = 0 and $n$ = 1.

-------------------------------

$C_{27}(1,3,8,10,17,19,24,26)$

$C_{27}(3,4,5,13,14,22,23,24)$

$C_{27}(2,3,7,11,16,20,24,25)$

-------------------------------

$T2_{54,3}(C_{54}(R^{54,1}_i))$, $p$ = 3, $x$ = 1, $y$ = 0 and $n$ = 2.

-------------------------------

$C_{54}(1,3,17,19,35,37,51,53)$

$C_{54}(3,7,11,25,29,43,47,51)$

$C_{54}(3,5,13,23,31,41,49,51)$

----------------------------------

\vspace{.2cm}
$T2_{54,3}(C_{54}(R^{54,2}_i))$, $p$ = 3, $x$ = 2, $y$ = 0 and $n$ = 2.

-------------------------------

$C_{54}(2,3,16,20,34,38,51,52)$

$C_{54}(3,4,14,22,32,40,50,51)$

$C_{54}(3,8,10,26,28,44,46,51)$

\vspace{.2cm}

xxxxxxxxxxxxxxxxxxxxxxxxxx

\vspace{.2cm}
$T2_{125,5}(C_{125}(R^{125,1}_i))$, $p$ = 5, $x$ = 1, $y$ = 0 and $n$ = 1.

-------------------------------

$C_{125}(1,5,24,26,49,51,74,76,99,101,120,124)$

$C_{125}(5,6,19,31,44,56,69,81,94,106,119,120)$

$C_{125}(5,11,14,36,39,61,64,86,89,111,114,120)$

$C_{125}(5,9,16,34,41,59,66,84,91,109,116,120)$

$C_{125}(4,5,21,29,46,54,71,79,96,104,120,121)$

----------------------------------

\vspace{.2cm}
$T2_{125,5}(C_{125}(R^{125,2}_i))$, $p$ = 5, $x$ = 2, $y$ = 0 and $n$ = 1.

-------------------------------

$C_{125}(2,5,23,27,48,52,73,77,98,102,120,123)$

$C_{125}(5,12,13,37,38,62,63,87,88,112,113,120)$

$C_{125}(3,5,22,28,47,53,72,78,97,103,120,122)$

$C_{125}(5,7,18,32,43,57,68,82,93,107,118,120)$

$C_{125}(5,8,17,33,42,58,67,83,92,108,117,120)$

\vspace{.2cm}

xxxxxxxxxxxxxxxxxxxxxxxxxx

\vspace{.2cm}
$T2_{250,5}(C_{250}(R^{250,1}_i))$, $p$ = 5, $x$ = 1, $y$ = 0 and $n$ = 2.

$C_{250}(1,5,49,51,99,101,149,151,199,201,245,249)$

$C_{250}(5,11,39,61,89,111,139,161,189,211,239,245)$

$C_{250}(5,21,29,71,79,121,129,171,179,221,229,245)$

$C_{250}(5,19,31,69,81,119,131,169,181,219,231,245)$

$C_{250}(5,9,41,59,91,109,141,159,191,209,241,245)$

----------------------------------

\vspace{.2cm}
$T2_{250,5}(C_{250}(R^{250,2}_i))$, $p$ = 5, $x$ = 2, $y$ = 0 and $n$ = 2.

----------------------------------

$C_{250}(2,5,48,52,98,102,148,152,198,202,245,248)$

$C_{250}(5,22,28,72,78,122,128,172,178,222,228,245)$

$C_{250}(5,8,42,58,92,108,142,158,192,208,242,245)$

$C_{250}(5,12,38,62,88,112,138,162,188,212,238,245)$

$C_{250}(5,18,32,68,82,118,132,168,182,218,232,245)$

----------------------------------

\vspace{.2cm}
$T2_{250,5}(C_{250}(R^{250,3}_i))$, $p$ = 5, $x$ = 3, $y$ = 0 and $n$ = 2.

----------------------------------

$C_{250}(3,5,47,53,97,103,147,153,197,203,245,247)$

$C_{250}(5,17,33,67,83,117,133,167,183,217,233,245)$

$C_{250}(5,13,37,63,87,113,137,163,187,213,237,245)$

$C_{250}(5,7,43,57,93,107,143,157,193,207,243,245)$

$C_{250}(5,23,27,73,77,123,127,173,177,223,227,245)$

----------------------------------

\vspace{.2cm}
$T2_{250,5}(C_{250}(R^{250,4}_i))$, $p$ = 5, $x$ = 4, $y$ = 0 and $n$ = 2.

----------------------------------

$C_{250}(4,5,46,54,96,104,146,154,196,204,245,246)$

$C_{250}(5,6,44,56,94,106,144,156,194,206,244,245)$

$C_{250}(5,16,34,66,84,116,134,166,184,216,234,245)$

$C_{250}(5,24,26,74,76,124,126,174,176,224,226,245)$

$C_{250}(5,14,36,64,86,114,136,164,186,214,236,245)$

\vspace{.2cm}

xxxxxxxxxxxxxxxxxxxxxxxxxxxxxxxxxxx

\vspace{.2cm}
$T2_{343,7}(C_{343}(R^{343,1}_i))$, $p$ = 7, $x$ = 1, $y$ = 0 and $n$ = 1.

---------------------------

$C_{343}(1,7,48,50,97,99,146,148,~195,197,244,246,293,295,336,342)$

$C_{343}(7,8,41,57,90,106,139,155,~188,204,237,253,286,302,335,336)$

$C_{343}(7,15,34,64,83,113,132,162,~181,211,230,260,279,309,328,336)$

$C_{343}(7,22,27,71,76,120,125,169,~174,218,223,267,272,316,321,336)$

$C_{343}(7,20,29,69,78,118,127,167,~176,216,225,265,274,314,323,336)$

$C_{343}(7,13,36,62,85,111,134,160,~183,209,232,258,281,307,330,336)$

$C_{343}(6,7,43,55,92,104,141,153,~190,202,239,251,288,300,336,337)$

----------------------------------------

\vspace{.2cm}
$T2_{343,7}(C_{343}(R^{343,2}_i))$, $p$ = 7, $x$ = 2, $y$ = 0 and $n$ = 1.

--------------------------

$C_{343}(2,7,47,51,96,100,145,149,~194,198,243,247,292,296,336,341)$

$C_{343}(7,16,33,65,82,114,131,163,~180,212,229,261,278,310,327,336)$

$C_{343}(7,19,30,68,79,117,128,166,~177,215,226,264,275,313,324,336)$

$C_{343}(5,7,44,54,93,103,142,152,~191,201,240,250,289,299,336,338)$

$C_{343}(7,9,40,58,89,107,138,156,~187,205,236,254,285,303,334,336)$

$C_{343}(7,23,26,72,75,121,124,170,~173,219,222,268,271,317,320,336)$

$C_{343}(7,12,37,61,86,110,135,159,~184,208,233,257,282,306,331,336)$

----------------------------------------

\vspace{.2cm}
$T2_{343,7}(C_{343}(R^{343,3}_i))$, $p$ = 7, $x$ = 3, $y$ = 0 and $n$ = 1.

-----------------------------------

$C_{343}(3,7,46,52,95,101,144,150,~193,199,242,248,291,297,336,340)$

$C_{343}(7,24,25,73,74,122,123,171,~172,220,221,269,270,318,319,336)$

$C_{343}(4,7,45,53,94,102,143,151,~192,200,241,249,290,298,336,339)$

$C_{343}(7,17,32,66,81,115,130,164,~ 179,213,228,262,277,311,326,336)$

$C_{343}(7,11,38,60,87,109,136,158,~185,207,234,256,283,305,332,336)$

$C_{343}(7,10,39,59,88,108,137,157,~186,206,235,255,284,304,333,336)$

$C_{343}(7,18,31,67,80,116,129,165,~178,214,227,263,276,312,325,336)$

\vspace{.2cm}

xxxxxxxxxxxxxxxxxxxxxxxxxx

\vspace{.2cm}
$T2_{686,7}(C_{686}(R^{686,1}_i))$, $p$ = 7, $x$ = 1, $y$ = 0 and $n$ = 2.

----------------------------

$C_{686}(1,7,97,99,195,197,293,295,~ 391,393,489,491,587,589,679,685)$

$C_{686}(7,15,83,113,181,211,279,309,~377,407,475,505,573,603,671,679)$

$C_{686}(7,29,69,127,167,225,265,323,~363,421,461,519,559,617,657,679)$

$C_{686}(7,43,55,141,153,239,251,337,~349,435,447,533,545,631,643,679)$

$C_{686}(7,41,57,139,155,237,253,335,~351,433,449,531,547,629,645,679)$

$C_{686}(7,27,71,125,169,223,267,321,~365,419,463,517,561,615,659,679)$

$C_{686}(7,13,85,111,183,209,281,307,~379,405,477,503,575,601,673,679)$

----------------------------------------

\vspace{.2cm}
$T2_{686,7}(C_{686}(R^{686,2}_i))$, $p$ = 7, $x$ = 2, $y$ = 0 and $n$ = 2.

--------------------------------------

$C_{686}(2,7,96,100,194,198,292,296,~390,394,488,492,586,590,679,684)$

$C_{686}(7,30,68,128,166,226,264,324,~362,422,460,520,558,618,656,679)$

$C_{686}(7,40,58,138,156,236,254,334,~352,432,450,530,548,628,646,679)$

$C_{686}(7,12,86,110,184,208,282,306,~380,404,478,502,576,600,674,679)$

$C_{686}(7,16,82,114,180,212,278,310,~376,408,474,506,572,604,670,679)$

$C_{686}(7,44,54,142,152,240,250,338,~348,436,446,534,544,632,642,679)$

$C_{686}(7,26,72,124,170,222,268,320,~366,418,464,516,562,614,660,679)$

----------------------------------------

\vspace{.2cm}
$T2_{686,7}(C_{686}(R^{686,3}_i))$, $p$ = 7, $x$ = 3, $y$ = 0 and $n$ = 2.

-----------------------------------

$C_{686}(3,7,95,101,193,199,291,297,~389,395,487,493,585,591,679,683)$

$C_{686}(7,45,53,143,151,241,249,339,~347,437,445,535,543,633,641,679)$

$C_{686}(7,11,87,109,185,207,283,305,~381,403,479,501,577,599,675,679)$

$C_{686}(7,31,67,129,165,227,263,325,~361,423,459,521,557,619,655,679)$

$C_{686}(7,25,73,123,171,221,269,319,~367,417,465,515,563,613,661,679)$

$C_{686}(7,17,81,115,179,213,277,311,~ 375,409,473,507,571,605,669,679)$

$C_{686}(7,39,59,137,157,235,255,333,~353,431,451,529,549,627,647,679)$

----------------------------------------

\vspace{.2cm}
$T2_{686,7}(C_{686}(R^{686,4}_i))$, $p$ = 7, $x$ = 4, $y$ = 0 and $n$ = 2.

---------------------------------

$C_{686}(4,7,94,102,192,200,290,298,~388,396,486,494,584,592,679,682)$

$C_{686}(7,38,60,136,158,234,256,332,~354,430,452,528,550,626,648,679)$

$C_{686}(7,18,80,116,178,214,276,312,~374,410,472,508,570,606,668,679)$

$C_{686}(7,24,74,122,172,220,270,318,~ 368,416,466,514,564,612,662,679)$

$C_{686}(7,32,66,130,164,228,262,326,~ 360,424,458,522,556,620,654,679)$

$C_{686}(7,10,88,108,186,206,284,304,~382,402,480,500,578,598,676,679)$

$C_{686}(7,46,52,144,150,242,248,340,~346,438,444,536,542,634,640,679)$

----------------------------------------

\vspace{.2cm}
$T2_{686,7}(C_{686}(R^{686,5}_i))$, $p$ = 7, $x$ = 5, $y$ = 0 and $n$ = 2.

--------------------------------------

$C_{686}(5,7,93,103,191,201,289,299,~387,397,485,495,583,593,679,681)$

$C_{686}(7,23,75,121,173,219,271,317,~369,415,467,513,565,611,663,679)$

$C_{686}(7,47,51,145,149,243,247,341,~ 345,439,443,537,541,635,639,679)$

$C_{686}(7,19,79,117,177,215,275,313,~373,411,471,509,569,607,667,679)$

$C_{686}(7,9,89,107,187,205,285,303,~383,401,481,499,579,597,677,679)$

$C_{686}(7,37,61,135,159,233,257,331,~355,429,453,527,551,625,649,679)$

$C_{686}(7,33,65,131,163,229,261,327,~359,425,457,523,555,621,653,679)$

----------------------------------------

\vspace{.2cm}
$T2_{686,7}(C_{686}(R^{686,6}_i))$, $p$ = 7, $x$ = 6, $y$ = 0 and $n$ = 2.

----------------------------------

$C_{686}(6,7,92,104,190,202,288,300,~ 386,398,484,496,582,594,679,680)$

$C_{686}(7,8,90,106,188,204,286,302,~ 384,400,482,498,580,596,678,679)$

$C_{686}(7,22,76,120,174,218,272,316,~370,414,468,512,566,610,664,679)$

$C_{686}(7,36,62,134,160,232,258,330,~356,428,454,526,552,624,650,679)$

$C_{686}(7,48,50,146,148,244,246,342,~344,440,442,538,540,636,638,679)$

$C_{686}(7,34,64,132,162,230,260,328,~358,426,456,524,554,622,652,679)$

$C_{686}(7,20,78,118,176,216,274,314,~ 372,412,470,510,568,608,666,679)$

\vspace{.2cm}

--------------xxxxx---------------

\end{document}